\documentclass[11pt]{article}
\usepackage{amsmath,amsfonts,amssymb,graphicx,epsfig,pdflscape,multirow,theorem} 
\usepackage{arydshln,cite}
\usepackage{simplewick}
\numberwithin{equation}{section}
\graphicspath{{./eps/}}
\long\def\ignore#1{}

\theorembodyfont{\itshape} 
\theoremheaderfont{\scshape}
\theoremstyle{plain}  
\newtheorem{Lemma}{Lemma}[section]

\newtheorem{Proposition}[Lemma]{Proposition}
\newtheorem{Conjecture}[Lemma]{Conjecture}
\newtheorem{Corollary}[Lemma]{Corollary}

\numberwithin{equation}{section}

\newcommand{\nc}{\newcommand}
\nc{\bib}{\bibitem}
\nc{\be}{\begin{equation}}
\nc{\ee}{\end{equation}}
\nc{\bea}{\begin{eqnarray}}
\nc{\eea}{\end{eqnarray}}
\nc{\nn}{\nonumber\\[.1cm]}
\nc{\chit}{\raisebox{0.25ex}{$\chi$}}
\nc{\chih}{\raisebox{0.25ex}{$\hat\chi$}}
\nc{\ch}{\mathrm{ch}}
\nc{\id}{\mathrm{id}}
\nc{\g}{\mathfrak{g}}
\nc{\Ac}{\mathcal{A}}
\nc{\Cc}{\mathcal{C}}
\nc{\Dc}{\mathcal{D}}
\nc{\Lc}{\mathcal{L}}
\nc{\Mc}{\mathcal{M}}
\nc{\Nc}{\mathcal{N}}
\nc{\Oc}{\mathcal{O}}
\nc{\Vc}{\mathcal{V}}
\nc{\Wc}{\mathcal{W}}
\nc{\pa}{\partial}
\nc{\al}{\alpha}
\nc{\eps}{\epsilon}
\nc{\Ga}{\Gamma}
\nc{\La}{\Lambda}
\nc{\la}{\lambda}
\nc{\om}{\omega}
\nc{\hc}{\mathfrak{h}}

\begin{document}

\topmargin -15mm
\oddsidemargin 05mm

\title{\mbox{}\vspace{0in}
\bf 
\huge
Layer structure\\[.1cm] 
of irreducible Lie algebra modules
\\[-.2cm]
}
\date{}
\author{}

\maketitle


\begin{center}
{\vspace{-8mm}\LARGE J{\o}rgen Rasmussen}
\\[.5cm]
{\em School of Mathematics and Physics, University of Queensland}\\
{\em St Lucia, Brisbane, Queensland 4072, Australia}
\\[.4cm] 
{\tt j.rasmussen\,@\,uq.edu.au}
\end{center}

\vspace{0.5cm}
\begin{abstract}
Let $\g$ be a finite-dimensional simple complex Lie algebra.
A layer sum is introduced as the sum of formal exponentials of the distinct weights appearing in an irreducible 
$\g$-module. It is argued that the character of every finite-dimensional irreducible $\g$-module
admits a decomposition in terms of layer sums, with only non-negative integer coefficients. 
Ensuing results include a new approach to the computation of Weyl characters and weight multiplicities,
and a closed-form expression for the number of distinct weights in a finite-dimensional irreducible $\g$-module. 
The latter is given by a polynomial in the Dynkin labels, of degree equal to the rank of $\g$.
\end{abstract}

\newpage
\tableofcontents

\newpage

\section{Introduction}
\label{Sec:Introduction}

Several classic results on the representation theory of Lie algebras~\cite{Hum72} are due to Weyl
and have been known for almost a century.
This includes his character and dimension formulas~\cite{Weyl} 
for finite-dimensional irreducible modules over simple complex Lie algebras of finite type.
These results are remarkably succinct and give fundamental insight into the structure 
of the modules. However, the character formula requires cumbersome manipulations to
reveal certain key details and does not offer a closed-form expression for the weight 
multiplicies~\cite{Fre54,Kos59}.
A primary objective of the present work is to find a new and computationally efficient way to obtain
descriptive expressions for the characters.

Let $\g$ be a finite-dimensional simple complex Lie algebra.
Instrumental to the approach presented here, it is asserted that the character of every finite-dimensional irreducible 
$\g$-module admits a decomposition in terms of so-called {\em layer sums}. Here, a layer sum 
is the sum of formal exponentials of the distinct weights appearing in an irreducible $\g$-module.
We find that the number of distinct weights is {\em polynomial}\, in the Dynkin labels of the highest weight 
characterising the finite-dimensional irreducible module, 
and that the degree of the polynomial is equal to the rank of $\g$.
Although some results on these numbers are known~\cite{Pos09},
their polynomial nature, in particular, does not appear to be discussed in the literature.

In a given finite-dimensional irreducible $\g$-module, the weight multiplicities are Weyl group invariant.
Determining its character thus amounts to specifying the multiplicities of the dominant integral weights 
appearing in the module, and working out the associated Weyl orbits. This is still a nontrivial task.
Here, it is proposed that the {\em inverse} problem, expressing the orbit sums in terms of
irreducible characters, has a simple solution. We thus assert that the orbit sum
corresponding to a dominant integral weight can be written as an alternating sum 
of finite-dimensional irreducible characters, where the sum is over the Weyl group.
Moreover, if the dominant integral weights are ordered according to their values under the layer polynomial,
these relations form an infinite linear system corresponding to a lower-triangular
matrix with $1$'s on the diagonal. For any $n\in\mathbb{N}$, one can then invert the top-left $n\times n$ part
of the matrix to obtain the `first' $n$ irreducible characters.

A similar approach can be applied to find explicit expressions for the layer sums in terms of irreducible characters.
In this case, the alternating sum is over an abelian group of order given by the number of
non-simple positive roots. Accordingly, a {\em reduced Weyl vector} appears in these expressions,
defined as half the sum of the non-simple positive roots.

In Section~\ref{Sec:Notation}, to fix our notation,
we review the basic Lie algebra theory needed in the subsequent sections.
We also introduce the notion of {\em auxiliary characters} to assist in the description of the relations
between orbit and layer sums and irreducible characters. This is based on a seemingly new polynomial identity
involving Weyl's dimension formula.

In Section~\ref{Sec:Decomp}, we discuss the {\em layer decomposition} of characters of finite-dimensional 
irreducible modules. Layer sums are introduced and their conjectured expressions in terms of irreducible 
characters are given. The corresponding layer polynomials are also defined and subsequently
expressed using Weyl's dimension formula. 
Some examples are presented, with additional ones deferred to Appendix~\ref{Sec:LayerEx}.
Layer polynomials can be constructed alternatively by counting the number of lattice points in 
the weight polytopes associated with the modules. That the two methods indeed agree is confirmed in 
Appendix~\ref{Sec:Equiv} for $A_2$, $A_3$, $B_2$, and $G_2$.

In Section~\ref{Sec:Characters}, we present the conjectured relations between orbit and layer sums 
and irreducible characters, including the weight multiplicities. 
As a corollary, we find that the order of the Weyl group can be written
as an alternating sum of the layer polynomial evaluated at points related by the shifted Weyl group action.
We use $G_2$ to illustrate the general results and to verify a nontrivial consistency condition.

Section~\ref{Sec:Discussion} contains some concluding remarks.

\section{Notation}
\label{Sec:Notation}

Let $X_r$ be a simple complex Lie algebra of finite type, where $X\in\{A,\ldots,G\}$ and 
$r=\mathrm{rank}\,X_r$. We denote the corresponding root system by $\Phi$, the set
of positive roots by $\Phi_+$, a base of simple roots by $\Delta=\{\al_1,\ldots,\al_r\}$, and the
set of non-simple positive roots by
\be
 \Phi_+':=\Phi_+\setminus\Delta.
\ee
The non-negative root lattice is defined as
\be
 Q_+:=\mathbb{N}_0\al_1+\ldots+\mathbb{N}_0\al_r,
\ee
while our convention for the Cartan matrix is as follows:
\be
 A=(A_{ij}),\qquad A_{ij}:=\langle\al_i^\vee,\al_j\rangle,\qquad i,j=1,\ldots,r.
\ee

Let $\hc=\mathrm{span}\{h_1,\ldots,h_r\}$ be a Cartan subalgebra of $X_r$. For $\la\in\hc^\ast$, we can write
\be
 \la=\la_1\om_1+\ldots+\la_r\om_r,
\ee
where $\{\om_1,\ldots,\om_r\}$ is the set of fundamental weights,
dual to the set of simple coroots, $\{\al_1^\vee,\ldots,\al_r^\vee\}$, 
while the scalars $\la_1,\ldots,\la_r$ are known as Dynkin labels.
Correspondingly, the respective sets of integral weights and of dominant integral weights are defined as
\be
 P:=\mathbb{Z}\,\om_1+\ldots+\mathbb{Z}\,\om_r,\qquad
 P_+:=\mathbb{N}_0\om_1+\ldots+\mathbb{N}_0\om_r.
\ee
The latter admits a partial ordering, where 
\be
 \mu\leq\la\quad\,\mathrm{if}\quad 
 \la-\mu\in Q_+.
\label{lamuA}
\ee

Finite-dimensional irreducible $X_r$-modules are exactly the irreducible highest-weight modules $L(\la)$ for which
$\la\in P_+$. For $\la,\mu\in P_+$, $\mu$ is a weight of $L(\la)$ if and only if $\mu\leq\la$.
For $\la\in P_+$, the set of distinct weights in $L(\la)$ is denoted by $P(\la)$, the set of distinct dominant integral
weights in $L(\la)$ by $P_+(\la)$, and the character of $L(\la)$ by $\ch_\la$, while
Weyl's dimension formula expresses the dimension of $L(\la)$ as
\be
 \dim L(\la)=\prod_{\al\in\Phi_+}\frac{\langle\al,\la+\rho\rangle}{\langle\al,\rho\rangle}.
\label{WeylDim}
\ee

Let $W$ denote the Weyl group associated with $X_r$, and $O_\la^{X_r}$, or simply $O_\la$, 
the corresponding Weyl orbit of $\la\in P$.
If $\mu\in P(\la)$, $\la\in P_+$, then so is every weight in $O_\mu$, and 
exactly one of the weights in $O_\mu$ is in $P_+$.
Simple Weyl reflections are denoted by $s_1,\ldots,s_r$, and $\ell(w)$ denotes the length of $w\in W$.
The shifted action of $w\in W$ on $\la$ is defined by
\be
 w\cdot\la:=w(\la+\rho)-\rho,
\ee
where
\be
 \rho:=\tfrac{1}{2}\sum_{\al\in\Phi_+}\!\al=\sum_{i=1}^r\om_i
\ee
is the Weyl vector. We introduce the {\em shifted Weyl orbit} of $\la\in P$ as
\be
 O_{\cdot\la}:=\{w\cdot\la\,|\,w\in W\}.
\ee

We shall be interested in the group
\be
 \mathbb{Z}_2^k\equiv(\mathbb{Z}_2)^{\times k},\qquad
 k=|\Phi_+'|=\tfrac{1}{2}(\dim X_r-3r),
\ee
and its `shifted action' on the zero weight, where, for each $\al\in\Phi_+'$,
the generator $z_\al\in\mathbb{Z}_2^k$ acts by subtracting the root $\al$:
\be
 z_\al\cdot0:=-\al.
\ee
For $\al\neq\al'$, $\al,\al'\in\Phi_+'$, the composition $z_\al z_{\al'}$ thus subtracts $\al+\al'$, whereas, by construction, $z_\al^2=\id$.
As for Weyl group elements, the length of $z\in\mathbb{Z}_2^k$ is denoted by $\ell(z)$ and defined as the number 
of basic generators (the ones of the form $z_\al$, $\al\in\Phi_+'$) appearing in a reduced decomposition of $z$.
Thus, the unique longest element, $\prod_{\al\in\Phi_+'}z_\al$, has length $k$ and acts by subtracting $2\rho'$, where 
the {\em reduced Weyl vector} $\rho'$ is defined as half the sum of the non-simple positive roots:
\be
 \rho':=\tfrac{1}{2}\sum_{\al\in\Phi_+'}\!\al
  =\sum_{i=1}^r(\om_i-\tfrac{1}{2}\al_i).
\ee

\subsection{Auxiliary characters}
\label{Sec:Auxiliary}

Expanded out, Weyl's dimension formula (\ref{WeylDim}) expresses $\dim L(\la)$ as a polynomial in the $r$ 
(non-negative integer) Dynkin labels. We denote by $D_{X_r}$, or simply $D$, the polynomial in the
$r$ variables $\la_1,\ldots,\la_r$ that agrees with $\dim L(\la)$ for $\la=(\la_1,\ldots,\la_r)\in P_+$.
By construction, $D$ is of degree $|\Phi_+|$. We will not distinguish between 
$D$ as a function of the $r$-tuple $(\la_1,\ldots,\la_r)$ and $D$ as a function of $\la\in\hc^\ast$, 
setting $D(\la_1,\ldots,\la_r)\equiv D(\la)$. For $G_2$, for example, the polynomial is given by
\be
 D_{G_2}(\la)
 =\tfrac{1}{120}(1+\la_1)(1+\la_2)(2+\la_1+\la_2)(3+2\la_1+\la_2)(4+3\la_1+\la_2)(5+3\la_1+2\la_2).
\ee
In general, $(1+\la_1)\ldots(1+\la_r)$ is a divisor of $D(\la)$.
It readily follows that $D(\la)=0$ if $\la_i=-1$ for some $i=1,\ldots,r$, for example.
\begin{Proposition}
For every $w\in W$ and $\la\in\hc^\ast$,
\be
 D(w\cdot\la)=(-1)^{\ell(w)}D(\la).
\ee
\end{Proposition}
\begin{Corollary}
Let $\la\in\hc^\ast$. If $D(\la)=0$, then $\la\not\in P_+$ and
\be
 \mu\in O_{\cdot\la}\quad\implies\quad D(\mu)=0.
\ee
\end{Corollary}
Following these results, for $\la\in P_+$, we introduce {\em auxiliary characters} for the
weights in $O_{\cdot\la}$ of the form $w\cdot\la\neq\la$, $w\in W$, as
\be
 \ch_{w\cdot\la}:=(-1)^{\ell(w)}\ch_{\la},
\ee
where it is noted that $w\cdot\la\neq\la$ implies $w\cdot\la\not\in P_+$.
In addition, if $\la\in P$ is a zero of $D$, then we set $\ch_\la:=0$.
For $G_2$,
\be
 s_1s_2\cdot(3\om_1-6\om_2)=\om_2,\qquad
 s_1s_2s_1\cdot(-4\om_1+4\om_2)=0,\qquad
 D_{G_2}(-3\om_1+\om_2)=0,
\ee
so we set
\be
 \ch_{3\om_1-6\om_2}=\ch_{\om_2},\qquad
 \ch_{-4\om_1+4\om_2}=-\ch_{0},\qquad
 \ch_{-3\om_1+\om_2}=0,
\ee
for example.
Although we do not provide details here, we note that the auxiliary characters can be understood using 
reflections about the edges of the fundamental Weyl chamber.

\section{Decomposition of irreducible modules}
\label{Sec:Decomp}

\subsection{Layer structure}
\label{Sec:Layer}

Let $\la\in P_+$. For each $\mu\in P(\la)$, let $v_\mu\in L(\la)$ be a vector in the $\hc$-eigenspace of eigenvalue 
$\mu$. That is, $v_\mu$ is a simultaneous eigenvector of the Cartan basis generators $\{h_1,\ldots,h_r\}$, 
with eigenvalues given by the Dynkin labels of $\mu$:
\be
 h_iv_\mu=\mu_iv_\mu,\qquad i=1,\ldots,r.
\ee
We refer to
\be
 \mathrm{span}\{v_\mu\,|\,\mu\in P(\la)\}
\label{Lcla}
\ee
as a {\em layer} corresponding to $\la$. A layer is thus a direct sum of one-dimensional $\hc$-modules, 
where the sum is over the weights in $P(\la)$. If $L(\la)$ contains $\hc$-eigenspaces of dimension greater 
than $1$, then the layer (\ref{Lcla}) is not an $X_r$-module nor unique. 
However, if an ordered basis is given for every one of the $\hc$-eigenspaces, we may consider
the unique layer formed by the first basis vectors.

Motivated by this, we introduce the {\em layer sum} $\Lc_\la^{X_r}$, or simply $\Lc_\la$,
as the sum of formal exponentials of the elements of $P(\la)$:
\be
 \Lc_\la:=\sum_{\mu\in P(\la)}e^\mu.
\ee
The next conjecture asserts that every irreducible character admits a {\em layer decomposition}
in terms of such layer sums.
\begin{Conjecture}
For $\la\in P_+$, 
\be
 \ch_\la=\sum_{\mu\in P_+(\la)}\!\!c_{\la,\mu}\Lc_\mu
\label{Lla}
\ee
for some $c_{\la,\mu}\in\mathbb{N}_0$.
\end{Conjecture}
Although any given layer sum $\Lc_\mu$ will appear in the decomposition of infinitely many distinct irreducible 
characters, it need not appear in the decomposition of $\ch_\la$ just because $\mu\in P_+(\la)$. 
In the case of $G_2$, for example, $0,\om_1\in P_+(\om_1+\om_2)$, but the layer decomposition
\be
 \ch_{\om_1+\om_2}=\Lc_{\om_1+\om_2}+\Lc_{2\om_2}+2\Lc_{\om_2}
\label{Lom1om2}
\ee
does not involve $\Lc_0$ nor $\Lc_{\om_1}$.
\begin{Corollary}
For $\la\in P_+$,
\be
 \dim L(\la)=\sum_{\mu\in P_+(\la)}\!\!c_{\la,\mu}|P(\mu)|.
\label{dimLasum}
\ee
\end{Corollary}
In the $G_2$ example above, we confirm that
\be
 \dim L(\om_1+\om_2)=|P(\om_1+\om_2)|+|P(2\om_2)|+2|P(\om_2)|=31+19+14=64.
\ee
\begin{Conjecture}
Up to permutations of summands, the layer decomposition (\ref{Lla}) is unique.
\end{Conjecture}
Thus, it is not only proposed that  layer sums play a fundamental role in the description
of finite-dimensional irreducible modules; they are in some sense canonical.

\subsection{Layer sums}
\label{Sec:LayerChar}

The following conjecture offers an explicit expression for the layer sums.
\begin{Conjecture}
\label{Con33}
For $\la\in P_+$,
\be
 \Lc_\la=\sum_{z\in\mathbb{Z}_2^k}(-1)^{\ell(z)}\ch_{\la+z\cdot0},
\label{chiLcla}
\ee
where summands $\ch_{\la+z\cdot0}$, for which $(\la+z\cdot0)\not\in P_+$, are interpreted as auxiliary characters.
\end{Conjecture}
To illustrate this conjecture, let us consider $G_2$. Our labelling convention is
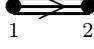
\setlength{\unitlength}{1mm}
\begin{center}
\begin{picture}(18,5)(-4.8,-2)
\thicklines
\put(0,0){\circle*{2}}
\put(10,0){\circle*{2}}
\put(0,-1){\line(1,0){10}}
\put(0,0){\line(1,0){10}}
\put(0,1){\line(1,0){10}}
\put(2.8,-1.5){{\LARGE $>$}}
\put(-0.7,-3.5){$_1$}
\put(9.1,-3.5){$_2$}
\end{picture}
\end{center}
in which case
\be
 A=\begin{pmatrix} 2&-1\\ -3&2\end{pmatrix}
\ee
and
\be
 \Phi_+'=\{\al_1+\al_2, \al_1+2\al_2, \al_1+3\al_2, 2\al_1+3\al_2\},\qquad
 \rho'=\tfrac{5}{2}\al_1+\tfrac{9}{2}\al_2=\tfrac{1}{2}\om_1+\tfrac{3}{2}\om_2.
\ee
The conjecture then asserts that
\begin{align}
 \Lc_\la^{G_2}&=\ch_{\la}
  -\big(\ch_{\la-\al_1-\al_2}+\ch_{\la-\al_1-2\al_2}+\ch_{\la-\al_1-3\al_2}+\ch_{\la-2\al_1-3\al_2}\big)
 \nn
 &+\big(\ch_{\la-2\al_1-3\al_2}+\ch_{\la-2\al_1-4\al_2}+\ch_{\la-3\al_1-4\al_2}+\ch_{\la-2\al_1-5\al_2}+\ch_{\la-3\al_1-5\al_2}
  +\ch_{\la-3\al_1-6\al_2}\big)
 \nn
 &-\big(\ch_{\la-3\al_1-6\al_2}+\ch_{\la-4\al_1-6\al_2}+\ch_{\la-4\al_1-7\al_2}+\ch_{\la-4\al_1-8\al_2}\big)
  +\ch_{\la-5\al_1-9\al_2}
 \nn
 &=(\ch_{\la}+\ch_{\la-5\al_1-9\al_2})
 -(\ch_{\la-\al_1-\al_2}+\ch_{\la-4\al_1-8\al_2})
 -(\ch_{\la-\al_1-2\al_2}+\ch_{\la-4\al_1-7\al_2})
 \nn
 &-(\ch_{\la-\al_1-3\al_2}+\ch_{\la-4\al_1-6\al_2})
 +(\ch_{\la-2\al_1-4\al_2}+\ch_{\la-3\al_1-5\al_2})
 +(\ch_{\la-3\al_1-4\al_2}+\ch_{\la-2\al_1-5\al_2}).
\label{LG2}
\end{align}
The rewriting in (\ref{LG2}) is due to simple cancellations of terms, and
indicates how the weights $\la-\mu$ and $\la-(2\rho'-\mu)$ can be paired up.
Similar rewritings are possible for all $X_r$, where the relative sign between $\ch_{\la-\mu}$
and $\ch_{\la-(2\rho'-\mu)}$ is given by the signature of the longest element of $\mathbb{Z}_2^k$. 
Since the length of that element equals $k=|\Phi_+'|$, the relative sign is given by $(-1)^{|\Phi_+'|}$.
In accordance with (\ref{LG2}), the relative sign for $G_2$ is $+1$.
For $A_2$, on the other hand, the relative sign is $-1$. Indeed, the number of non-simple positive roots for 
$A_2$ is $|\Phi_+'|=1$, the reduced Weyl vector is given by
$\rho'=\frac{1}{2}(\al_1+\al_2)$, and the layer sums are given by
\be
 \Lc_\la^{A_2}=\ch_\la-\ch_{\la-\al_1-\al_2}.
\label{LA2}
\ee

\subsection{Layer polynomial}
\label{Sec:LayerPol}

For $\la\in P_+$, the number of distinct weights in $L(\la)$
is given by $|P(\la)|$. This may be computed as a weighted sum over the elements in
$P_+(\la)$, weighting the elements by the corresponding orbit lengths, as
\be
 |P(\la)|=\sum_{\mu\in P_+(\la)}\!\!|O_\mu|.
\label{Ga}
\ee
The expression (\ref{chiLcla}) implies the following alternative expression for $|P(\la)|$.
\begin{Corollary}
For $\la\in P_+$,
\be
 |P(\la)|=\sum_{z\in\mathbb{Z}_2^k}(-1)^{\ell(z)}D(\la+z\cdot0).
\label{dimLc}
\ee
\end{Corollary}
This is a polynomial in the $r$ (non-negative integer) Dynkin labels.
We denote by $R_{X_r}$, or simply $R$, the polynomial in the
$r$ variables $\la_1,\ldots,\la_r$ that agrees with the expression in (\ref{dimLc}) for all $\la=(\la_1,\ldots,\la_r)\in P_+$, 
and refer to it as the corresponding {\em layer polynomial}. That is,
\be
 R(\la):=\sum_{z\in\mathbb{Z}_2^k}(-1)^{\ell(z)}D(\la+z\cdot0),\qquad \la\in\hc^\ast.
\ee
As indicated, we are not distinguishing between $R$ as a function of the $r$-tuple $(\la_1,\ldots,\la_r)$ and $R$ 
as a function of $\la\in\hc^\ast$, setting $R(\la_1,\ldots,\la_r)\equiv R(\la)$.
In Appendix~\ref{Sec:Equiv}, we verify that (\ref{Ga}) and (\ref{dimLc}) agree for $A_2$, $A_3$, $B_2$, and $G_2$,
thereby providing evidence for Conjecture~\ref{Con33}. By construction, the coefficients in $R(\la)$ are all rational.
\begin{Conjecture}
The polynomial $R(\la)$ has degree $r$, contains {\scriptsize $\begin{pmatrix} 2r\\ r\end{pmatrix}$} 
distinct terms, and has only positive coefficients.
\end{Conjecture}
As the maximum number of distinct terms in a polynomial of degree $n$ in $m$ variables is 
{\scriptsize $\begin{pmatrix} n+m\\ m\end{pmatrix}$}, it is thus asserted that this bound is saturated for all 
$R_{X_r}(\la)$. For $A_2$ and $G_2$, for example, the layer polynomials are found to be given by
\be
 R_{A_2}(\la)=1
  +\frac{3}{2}(\la_1+\la_2)
  +\frac{1}{2}(\la_1^2+\la_2^2+4\la_1\la_2)
\label{PA2}
\ee
and
\be
 R_{G_2}(\la)=1+3(\la_1+\la_2)+3(3\la_1^2+\la_2^2+4\la_1\la_2),
\label{PG2}
\ee
both having degree $2$ and containing $6$ distinct terms with only positive coefficients.
As an aside, we can rewrite $R_{G_2}$ as
\be
 R_{G_2}(\la)=1+3(\la_1+\la_2)(1+3\la_1+\la_2)
 =1+6\big[\!\begin{pmatrix} \la_1+\la_2+1\\ 2\end{pmatrix} +\la_1(\la_1+\la_2)\big],
\label{PG2mod6}
\ee
showing that $R_{G_2}$ evaluated at integer arguments gives $1$ plus an integer multiple of $6$.
Moreover, if one or more of the Dynkin labels is $0$, the polynomial expression $R(\la)$ simplifies considerably.
Particularly compact such specialised polynomials are
\be
 R_{A_r}(n\om_1)=R_{A_r}(n\om_r)=\begin{pmatrix} r+n\\ r\end{pmatrix},\qquad
  R_{B_r}(n\om_r)=(1+n)^r.
\ee

\section{Character expressions}
\label{Sec:Characters}

Here, we present two new ways of computing characters of finite-dimensional irreducible $X_r$-modules,
and a new expression for the weight multiplicities.

\subsection{Orbit sums}
\label{Sec:OrbitChar}

For $\la\in P_+$, the orbit sum $m_\la$ is defined as the sum of formal exponentials of the elements of $O_\la$, 
that is,
\be
 m_\la:=\sum_{\mu\in O_\la}e^\mu.
\label{chOla}
\ee
To specify the Lie algebra, we may write $m_\la^{X_r}$.
Essentially by construction, a layer sum can be expressed in terms of orbit sums as
\be
 \Lc_\la=\sum_{\mu\in P_+(\la)}\!\!m_\mu,
\label{LcO}
\ee
with all multiplicities being $1$.
\begin{Conjecture}
\label{Con41}
For $\la\in P_+$,
\be
 m_\la=\frac{|O_\la|}{|W|}\sum_{w\in W}(-1)^{\ell(w)}\ch_{\la+w\cdot0},
\label{chO}
\ee
where summands $\ch_{\la+w\cdot0}$, for which $(\la+w\cdot0)\not\in P_+$, are interpreted as auxiliary characters.
\end{Conjecture}
For $A_2$ and $G_2$, for example, we find that
\begin{align}
 m_\la^{A_2}
 &=\frac{|O_\la|}{6}\big[\ch_{\la}
   -(\ch_{\la-\al_1}+\ch_{\la-\al_2})
  +(\ch_{\la-2\al_1-\al_2}+\ch_{\la-\al_1-2\al_2})
  -\ch_{\la-2\al_1-2\al_2}\big]
 \nn
 &=\frac{|O_\la|}{6}\big[(\ch_{\la}-\ch_{\la-2\al_1-2\al_2})-(\ch_{\la-\al_1}-\ch_{\la-\al_1-2\al_2})
  -(\ch_{\la-\al_2}-\ch_{\la-2\al_1-\al_2})\big]
\end{align}
and
\begin{align}
 m_\la^{G_2}
 &=\frac{|O_\la|}{12}\big[\ch_{\la}
   -(\ch_{\la-\al_1}+\ch_{\la-\al_2})
  +(\ch_{\la-\al_1-4\al_2}+\ch_{\la-2\al_1-\al_2})
  -(\ch_{\la-4\al_1-4\al_2}+\ch_{\la-2\al_1-6\al_2})
 \nn
 &+(\ch_{\la-4\al_1-9\al_2}+\ch_{\la-5\al_1-6\al_2})
  -(\ch_{\la-6\al_1-9\al_2}+\ch_{\la-5\al_1-10\al_2})+\ch_{\la-6\al_1-10\al_2}\big]
 \nn
 &=\frac{|O_\la|}{12}\big[(\ch_{\la}+\ch_{\la-6\al_1-10\al_2})-(\ch_{\la-\al_1}+\ch_{\la-5\al_1-10\al_2})
  -(\ch_{\la-\al_2}+\ch_{\la-6\al_1-9\al_2})
 \nn
 &
  +(\ch_{\la-2\al_1-\al_2}+\ch_{\la-4\al_1-9\al_2})
  +(\ch_{\la-\al_1-4\al_2}+\ch_{\la-5\al_1-6\al_2})
  -(\ch_{\la-4\al_1-4\al_2}+\ch_{\la-2\al_1-6\al_2})\big].
\end{align}
The rewritings indicate how the weights $\la-\mu$ and $\la-(2\rho-\mu)$ can be paired up.
Similar rewritings are possible for all $X_r$, where the relative sign between $\ch_{\la-\mu}$
and $\ch_{\la-(2\rho-\mu)}$ is given by the signature of the longest element of $W$. Since the length
of that element equals $|\Phi_+|$, the relative sign is given by $(-1)^{|\Phi_+|}$.
As an illustration of how simplifications may be possible, for $B_2$, Conjecture~\ref{Con41} asserts that
\begin{align}
 m_{2\om_2}^{B_2}
 &=\tfrac{1}{2}\big[
  \ch_{2\om_2}-\ch_{-2\om_1+4\om_2}-\ch_{\om_1}+\ch_{\om_1-2\om_2}+\ch_{-3\om_1+4\om_2}
  -\ch_{-3\om_1+2\om_2}-\ch_{-2\om_2}+\ch_{-2\om_1}\big]
 \nn
 &=\tfrac{1}{2}\big[
  \ch_{2\om_2}-(-\ch_{2\om_2})-\ch_{\om_1}+(-\ch_{0})+(-\ch_{\om_1})
  -\ch_{0}-0+0\big]
 \nn
 &=\ch_{2\om_2}-\ch_{\om_1}-\ch_{0}.
\label{m2w2B2}
\end{align}
This is seen to agree with the orbit sum
\be
  m_{2\om_2}^{B_2}
=e^{2\om_2}+e^{-2\om_1+2\om_2}+e^{2\om_1-2\om_2}+e^{-2\om_2}
\ee
computed using (\ref{chOla}).

The expression (\ref{chO}) implies the following polynomial identity.
\begin{Corollary}
For $\la\in\hc^\ast$,
\be
 \sum_{w\in W}(-1)^{\ell(w)}D(\la+w\cdot0)=|W|.
\label{sumW}
\ee
\end{Corollary}
Despite its appearance, the sum in (\ref{sumW}) is thus found to be {\em independent} of $\la$.

\subsection{Irreducible characters as sums of orbit sums}
\label{Sec:IrrCharOrbit}

For $\la\in P_+$, the character of $L(\la)$ is of the form
\be
 \ch_\la=\sum_{\mu\in P}m_{\la,\mu}e^\mu
 =\sum_{\mu\in P_+} \!m_{\la,\mu} m_\mu,
\label{chla}
\ee
where the {\em weight multiplicities} $m_{\la,\mu}$ are non-negative integers. 
It is well known that $m_{\la,\mu}=0$ unless 
$\mu\leq\la$, but, for later convenience, we let the summation in (\ref{chla}) be over $\mu\in P_+$. 
Weyl's character formula expresses the character as
\be
 \ch_\la=\frac{\sum_{w\in W}(-1)^{\ell(w)}e^{w\cdot\la}}{\prod_{\al\in\Phi_+}(1-e^{-\al})}.
\ee

As discussed in the following, we find that the expressions (\ref{chO}) for orbit sums can be inverted.
This yields a straightforward approach to the computation of irreducible characters, including the weight 
multiplicities. First, we say that a pair of weights $\la,\mu\in P_+$ are related as follows:
\be
 \mu\prec\la\quad\mathrm{if}\quad R(\mu)<R(\la),\qquad\quad
 \mu\preceq\la\quad\mathrm{if}\quad R(\mu)\leq R(\la).
\ee
\begin{Proposition}
Let $\la,\mu\in P_+$. Then,
\be
 \mu<\la\ \implies\ \mu\prec\la.
\ee
\end{Proposition}
Second, choose an ordering $\Oc$ of the elements of $P_+$ such that $\mu$ appears before $\la$ if 
$\mu\preceq\la$. This is always possible, although the ordering need not be unique. For example, since
$R_{G_2}(3\om_1)=R_{G_2}(5\om_2)$, the ordering is not unique in the case of $G_2$.
Third, let 
\be
 \Mc=(m_{\la,\mu}),\qquad \la,\mu\in\Oc,
\ee
denote the infinite-dimensional matrix whose entries are given by the multiplicities in the last expression 
in (\ref{chla}), with the weights labelling the rows and columns ordered as in $\Oc$.
\begin{Conjecture}
$\Mc$ is a lower-triangular matrix with $1$'s on the diagonal.
\end{Conjecture}
\begin{Corollary}
The row of $\Mc^{-1}$ that corresponds to $\la\in P_+$ is read off (\ref{chO}).
\end{Corollary}
As a consequence, the {\em entire family} of irreducible characters $\ch_\la$, $\la\in P_+$, is obtained by
inverting the matrix $\Mc^{-1}$. 
Due to the triangular structure of $\Mc^{-1}$, we may choose to compute the finite set 
$\{\ch_\mu\,|\,\mu\preceq\la\}$ for any given $\la\in P_+$. 
This involves `finitising' $\Oc$ to include only the terms up to and including $\la$. We denote the ensuing 
ordered set by 
$\Oc_\la$. The corresponding top-left $|\Oc_\la|\times|\Oc_\la|$ part of $\Mc^{-1}$ is denoted by $\Mc^{-1}_\la$.

In this regard, we note that the inverse of a lower-triangular matrix $B=(b_{ij})$ with $1$'s on the diagonal
is a matrix of the same type, with
\be
 (B^{-1})_{ij}=-b_{ij}-\sum_{k=1}^{i-j-1}(-1)^k\sum_{j<\ell_k<\ell_{k-1}<\ldots<\ell_1<i}
  b_{i\ell_1}b_{\ell_1\ell_2}\ldots b_{\ell_kj},\qquad i>j.
\ee
It follows, in particular, that, if the entries of $B$ are all integer, then so are the entries of $B^{-1}$. 
Despite the division by $|W|$ in (\ref{chO}), the coefficient to any $\ch_\mu$ in the final expression
is therefore integer, as illustrated in (\ref{m2w2B2}).

As an example, let us consider $G_2$ and focus on the computation of 
$\{\ch_\mu\,|\,\mu\preceq2\om_1+2\om_2\}$. The corresponding `finitisation' of $\Oc$ is given by
\be
 \Oc_{2\om_1+2\om_2}=\{0,\om_2,\om_1,2\om_2,\om_1+\om_2,3\om_2,2\om_1,\om_1+2\om_2,4\om_2,
  2\om_1+\om_2,\om_1+3\om_2,3\om_1,5\om_2,2\om_1+2\om_2\},
\label{O2w12w2}
\ee
where the only freedom was the choice to place $3\om_1$ before $5\om_2$. Using (\ref{chO}), we find
\be
 \Mc^{-1}_{2\om_1+2\om_2}=\left(\!\!\begin{array}{rrrrrrrrrrrrrr}
 1&0&0&0&0&0&0&0&0&0&0&0&0&0\\
 -1&1&0&0&0&0&0&0&0&0&0&0&0&0\\
 -1&-1&1&0&0&0&0&0&0&0&0&0&0&0\\
 0&-1&-1&1&0&0&0&0&0&0&0&0&0&0\\
 2&0&0&-2&1&0&0&0&0&0&0&0&0&0\\
 0&1&-1&0&-1&1&0&0&0&0&0&0&0&0\\
 -1&1&0&0&0&-1&1&0&0&0&0&0&0&0\\
 0&0&1&1&-1&-1&-1&1&0&0&0&0&0&0\\
 -1&0&1&0&0&0&0&-1&1&0&0&0&0&0\\
 0&-1&0&1&0&1&0&-1&-1&1&0&0&0&0\\
 2&-1&-1&0&2&0&-1&0&-1&-1&1&0&0&0\\
 -1&0&0&0&0&0&0&0&1&0&-1&1&0&0\\
 -1&0&0&0&0&0&1&0&0&0&-1&0&1&0\\
 1&1&0&-2&0&1&0&1&0&-1&0&-1&-1&1
 \end{array}\!\right),
\ee
where the $0$ in position $(13,12)$ confirms the freedom to re-order $3\om_1$ and $5\om_2$.
Inverting the matrix yields the character expressions
\begin{align}
 \ch_{0}&=m_0
 \nn
 \ch_{\om_2}&=m_{\om_2}+m_{0}
 \nn
 \ch_{\om_1}&=m_{\om_1}+m_{\om_2} +2m_{0}
 \nn
 \ch_{2\om_2}&=m_{2\om_2}+m_{\om_1}+2m_{\om_2} +3m_{0}
 \nn
 \ch_{\om_1+\om_2}&=m_{\om_1+\om_2}+2m_{2\om_2}+2m_{\om_1}
  +4m_{\om_2} +4m_{0}
 \nn
 \ch_{3\om_2}&=m_{3\om_2}+m_{\om_1+\om_2}+2m_{2\om_2}+3m_{\om_1}
  +4m_{\om_2} +5m_{0}
 \nn
 \ch_{2\om_1}&=m_{2\om_1}+m_{3\om_2}+m_{\om_1+\om_2}+2m_{2\om_2}
  +3m_{\om_1}+3m_{\om_2} +5m_{0}
 \nn
 \ch_{\om_1+2\om_2}&=m_{\om_1+2\om_2}+m_{2\om_1}+2m_{3\om_2}+3m_{\om_1+\om_2}
  +5m_{2\om_2}+6m_{\om_1}+8m_{\om_2} +9m_{0}
 \nn
 \ch_{4\om_2}&=m_{4\om_2}+m_{\om_1+2\om_2}+m_{2\om_1}+2m_{3\om_2}
  +3m_{\om_1+\om_2}+5m_{2\om_2}+5m_{\om_1}+7m_{\om_2}+8m_{0}
 \nn
 \ch_{2\om_1+\om_2}&=m_{2\om_1+\om_2}+m_{4\om_2}+2m_{\om_1+2\om_2}+2m_{2\om_1}
  +3m_{3\om_2}+5m_{\om_1+\om_2}+7m_{2\om_2}+7m_{\om_1}+10m_{\om_2}
 \nn
 &+10m_{0}
 \nn
 \ch_{\om_1+3\om_2}&=m_{\om_1+3\om_2}+m_{2\om_1+\om_2}+2m_{4\om_2}
  +3m_{\om_1+2\om_2}+4m_{2\om_1}+6m_{3\om_2}+7m_{\om_1+\om_2}+10m_{2\om_2}
 \nn
 &+12m_{\om_1}+14m_{\om_2}+16m_{0}
 \nn
 \ch_{3\om_1}&=m_{3\om_1}+m_{\om_1+3\om_2}+m_{2\om_1+\om_2}+m_{4\om_2}
  +2m_{\om_1+2\om_2}+3m_{2\om_1}+4m_{3\om_2}+4m_{\om_1+\om_2}+5m_{2\om_2}
 \nn
 &+7m_{\om_1}+7m_{\om_2}+9m_{0}
 \nn
 \ch_{5\om_2}&=m_{5\om_2}+m_{\om_1+3\om_2}+m_{2\om_1+\om_2}+2m_{4\om_2}
  +3m_{\om_1+2\om_2}+3m_{2\om_1}+5m_{3\om_2}+6m_{\om_1+\om_2}+8m_{2\om_2}
 \nn
 &+9m_{\om_1}+11m_{\om_2}+12m_{0}
 \nn
 \ch_{2\om_1+2\om_2}&=m_{2\om_1+2\om_2}+m_{5\om_2}+m_{3\om_1}
  +2m_{\om_1+3\om_2}+3m_{2\om_1+\om_2}+4m_{4\om_2}
  +6m_{\om_1+2\om_2}+7m_{2\om_1}+9m_{3\om_2}
 \nn
 &+11m_{\om_1+\om_2}+15m_{2\om_2}
  +16m_{\om_1}+19m_{\om_2}+21m_{0}.
\label{chmG2}
\end{align}
As the orbit sums are readily worked out, we have thus obtained a whole family of irreducible characters
by computing the inverse of a simple, integer, lower-triangular matrix with $1$'s on the diagonal.

\subsection{Irreducible characters as sums of layer sums}
\label{Sec:IrrCharLayer}

We find that the expressions (\ref{chiLcla}) for layer sums can be inverted, allowing us to write an
irreducible character as a sum of layer sums.
As in Section~\ref{Sec:IrrCharOrbit}, choose an ordering $\Oc$ of the elements of $P_+$, and
let 
\be
 \Cc=(c_{\la,\mu}),\qquad \la,\mu\in\Oc,
\ee
denote the infinite-dimensional matrix whose entries are given by the multiplicities in
the decomposition
\be
 \ch_\la=\sum_{\mu\in P_+}\!c_{\la,\mu}\Lc_\mu.
\ee
According to (\ref{Lla}), the summation could be restricted to $\mu\in P_+(\la)$, but it is convenient to let it be 
over all of $P_+$, with $c_{\la,\mu}=0$ if $\mu\not\in P_+(\la)$.
\begin{Conjecture}
$\Cc$ is a lower-triangular matrix with $1$'s on the diagonal.
\end{Conjecture}
\begin{Corollary}
The row of $\Cc^{-1}$ that corresponds to $\la\in P_+$ is read off (\ref{chiLcla}).
\end{Corollary}
As a consequence, the entire family of irreducible characters $\ch_\la$, $\la\in P_+$, is obtained 
as expressions in layer sums by inverting the matrix $\Cc^{-1}$. 
As before, due to the triangular structure of $\Cc^{-1}$, we may choose to compute the finite set 
$\{\ch_\mu\,|\,\mu\preceq\la\}$ for any given $\la\in P_+$. 
The corresponding top-left $|\Oc_\la|\times|\Oc_\la|$ part of $\Cc^{-1}$ is denoted by $\Cc^{-1}_\la$.

To illustrate, let us again consider the computation of $\{\ch_\mu\,|\,\mu\preceq2\om_1+2\om_2\}$ for $G_2$. 
Using (\ref{chiLcla}), relative to the ordered set $\Oc_{2\om_1+2\om_2}$ given in (\ref{O2w12w2}), we find
\be
 \Cc^{-1}_{2\om_1+2\om_2}=\left(\!\!\begin{array}{rrrrrrrrrrrrrr}
 1&0&0&0&0&0&0&0&0&0&0&0&0&0\\
 0&1&0&0&0&0&0&0&0&0&0&0&0&0\\
 -1&0&1&0&0&0&0&0&0&0&0&0&0&0\\
 -1&-1&0&1&0&0&0&0&0&0&0&0&0&0\\
 1&-1&0&-1&1&0&0&0&0&0&0&0&0&0\\
 1&0&-1&-1&0&1&0&0&0&0&0&0&0&0\\
 0&1&-1&-1&0&0&1&0&0&0&0&0&0&0\\
 0&1&0&0&-1&-1&0&1&0&0&0&0&0&0\\
 -1&1&1&0&-1&-1&0&0&1&0&0&0&0&0\\
 -1&0&1&1&-1&0&0&-1&0&1&0&0&0&0\\
 1&-1&0&1&1&0&-1&-1&-1&0&1&0&0&0\\
 0&-1&0&1&1&0&-1&-1&0&0&0&1&0&0\\
 0&-1&0&1&1&0&0&-1&-1&0&0&0&1&0\\
 0&0&0&-1&1&1&0&0&0&-1&-1&0&0&1
 \end{array}\!\right).
\ee
Inverting this matrix yields the layer decompositions
\begin{align}
 \ch_0&=\Lc_0
 \nn
 \ch_{\om_2}&=\Lc_{\om_2}
 \nn
 \ch_{\om_1}&=\Lc_{\om_1}+\Lc_0
 \nn
 \ch_{2\om_2}&=\Lc_{2\om_2}+\Lc_{\om_2}+\Lc_0
 \nn
 \ch_{\om_1+\om_2}&=\Lc_{\om_1+\om_2}+\Lc_{2\om_2}+2\Lc_{\om_2}
 \nn
 \ch_{3\om_2}&=\Lc_{3\om_2}+\Lc_{2\om_2}+\Lc_{\om_1}+\Lc_{\om_2}+\Lc_0
 \nn
 \ch_{2\om_1}&=\Lc_{2\om_1}+\Lc_{2\om_2}+\Lc_{\om_1}+2\Lc_0
 \nn
 \ch_{\om_1+2\om_2}&=\Lc_{\om_1+2\om_2}+\Lc_{3\om_2}+\Lc_{\om_1+\om_2}
  +2\Lc_{2\om_2}+\Lc_{\om_1}+2\Lc_{\om_2}+\Lc_0
 \nn
 \ch_{4\om_2}&=\Lc_{4\om_2}+\Lc_{3\om_2}+\Lc_{\om_1+\om_2}+2\Lc_{2\om_2}
  +2\Lc_{\om_2}+\Lc_0
 \nn
 \ch_{2\om_1+\om_2}&=\Lc_{2\om_1+\om_2}+\Lc_{\om_1+2\om_2}+\Lc_{3\om_2}
  +2\Lc_{\om_1+\om_2}+2\Lc_{2\om_2}+3\Lc_{\om_2}
 \nn
 \ch_{\om_1+3\om_2}&=\Lc_{\om_1+3\om_2}+\Lc_{4\om_2}+\Lc_{\om_1+2\om_2}
  +\Lc_{2\om_1}+2\Lc_{3\om_2}+\Lc_{\om_1+\om_2}+3\Lc_{2\om_2}
  +2\Lc_{\om_1}+2\Lc_{\om_2}+2\Lc_0
 \nn
 \ch_{3\om_1}&=\Lc_{3\om_1}+\Lc_{\om_1+2\om_2}+\Lc_{2\om_1}+\Lc_{3\om_2}
  +\Lc_{2\om_2}+2\Lc_{\om_1}+2\Lc_0
 \nn
 \ch_{5\om_2}&=\Lc_{5\om_2}+\Lc_{4\om_2}+\Lc_{\om_1+2\om_2}+2\Lc_{3\om_2}
  +\Lc_{\om_1+\om_2}+2\Lc_{2\om_2}+\Lc_{\om_1}+2\Lc_{\om_2}
  +\Lc_0
 \nn
 \ch_{2\om_1+2\om_2}&=\Lc_{2\om_1+2\om_2}+\Lc_{\om_1+3\om_2}+\Lc_{2\om_1+\om_2}
  +\Lc_{4\om_2}+2\Lc_{\om_1+2\om_2}+\Lc_{2\om_1}+2\Lc_{3\om_2}
  +2\Lc_{\om_1+\om_2}
 \nn
 &+4\Lc_{2\om_2}+\Lc_{\om_1}+3\Lc_{\om_2}+2\Lc_0.
\label{chLG2}
\end{align}

\subsection{Weight multiplicities}
\label{Sec:DomM}

The relation (\ref{LcO}) is readily extended from a sum over $P_+(\la)$ to a sum over all of $P_+$.
For $\la\in P_+$, we may thus write
\be
 \Lc_\la=\sum_{\mu\in P_+}\Dc_{\la,\mu}m_\mu,
\ee
where the {\em dominance matrix} $\Dc$ has entries
\be
 \Dc_{\la,\mu}=\left\{\begin{array}{c} \!\!1,\quad\, \mu\in P_+(\la), \\[.25cm] 
  \!\!0,\quad\, \mu\not\in P_+(\la). \end{array}\right.
\label{Dc}
\ee
Relative to an ordering of $P_+$ of the form $\Oc$ discussed in Section~\ref{Sec:IrrCharOrbit},
this is clearly a lower-triangular matrix with $1$'s on the diagonal.
Combining the conjectures above then implies the following relation.
\begin{Corollary}
\be
 \Cc^{-1}\Mc=\Dc.
\label{Cm1MD}
\ee
\end{Corollary}
It follows that, for $\la,\mu\in P_+$, the weight multiplicity $m_{\la,\mu}$ can be expressed as
\be
 m_{\la,\mu}=(\Cc\Dc)_{\la,\mu}=\sum_{\nu\in P_+}c_{\la,\nu}\Dc_{\nu,\mu}=\sum_{\mu\leq\nu\leq\la}c_{\la,\nu},
\label{m}
\ee
in accordance with the fact that $m_{\la,\mu}=0$ unless $\mu\leq\la$.

Viewing (\ref{Cm1MD}) as a consistency condition, let us verify it in the $G_2$ example above. We thus 
compute the $14\times14$ matrix product
\be
 \Cc_{2\om_1+2\om_2}^{-1}\Mc_{2\om_1+2\om_2}=\left(\!\begin{array}{rrrrrrrrrrrrrr}
 1&0&0&0&0&0&0&0&0&0&0&0&0&0\\
 1&1&0&0&0&0&0&0&0&0&0&0&0&0\\
 1&1&1&0&0&0&0&0&0&0&0&0&0&0\\
 1&1&1&1&0&0&0&0&0&0&0&0&0&0\\
 1&1&1&1&1&0&0&0&0&0&0&0&0&0\\
 1&1&1&1&1&1&0&0&0&0&0&0&0&0\\
 1&1&1&1&1&1&1&0&0&0&0&0&0&0\\
 1&1&1&1&1&1&1&1&0&0&0&0&0&0\\
 1&1&1&1&1&1&1&1&1&0&0&0&0&0\\
 1&1&1&1&1&1&1&1&1&1&0&0&0&0\\
 1&1&1&1&1&1&1&1&1&1&1&0&0&0\\
 1&1&1&1&1&1&1&1&1&1&1&1&0&0\\
 1&1&1&1&1&1&1&1&1&1&1&0&1&0\\
 1&1&1&1&1&1&1&1&1&1&1&1&1&1
 \end{array}\!\right).
\ee
Noting that the zero in position $(13,12)$ reflects that 
$3\om_1\not\in P_+(5\om_2)$, in accordance with (\ref{Dc}), this is indeed seen to confirm (\ref{Cm1MD}).
In the same example, one may verify the expression (\ref{m}) for the weight multiplicities. 
In particular, from (\ref{chLG2}), we find the partial row sums
\be
 m_{2\om_1+\om_2,\om_1}=\sum_{\om_1\leq\nu\leq2\om_1+\om_2}c_{2\om_1+\om_2,\nu}
  =0+2+2+1+0+1+0+1=7
\ee
and
\be
 m_{2\om_1+2\om_2,0}=\sum_{0\leq\nu\leq2\om_1+2\om_2}c_{2\om_1+2\om_2,\nu}
  =2+3+1+4+2+2+1+2+1+1+1+0+0+1=21,
\ee
in accordance with (\ref{chmG2}).

\section{Discussion}
\label{Sec:Discussion}

Layers have been introduced to describe finite-dimensional irreducible $X_r$-modules.
This has allowed us to devise new methods for computing Weyl characters and weight multiplicities, 
including whole families of characters at a time, and to find a polynomial giving the number of distinct weights 
in such an $X_r$-module. We also expect to be able to construct closed-form expressions for
the weight multiplicities, and that the layer structure will enable the determination of explicit bases for the modules.
We hope to return elsewhere with a discussion of these problems and with proofs of the various conjectures 
put forward in the present work. It seems natural to expect that the related and well-developed theory of symmetric 
functions~\cite{Mac95} may play a role in such proofs, at least for the $A$-series.
We also intend to study how our new insight and results extend to infinite-dimensional
modules and to the representation theory of Lie superalgebras and affine Lie algebras.

\subsection*{Acknowledgements}

This work was supported by the Australian Research Council under the Discovery Project scheme, 
project number DP160101376. The author is grateful to the mathematical research institute MATRIX in Australia 
and the Centre for the Mathematics of Quantum Theory (QMATH) at the University of Copenhagen, 
where parts of the work were carried out.
The author thanks Mark Gould, Jesper Grodal, Phil Isaac, Hans Plesner Jakobsen, Masoud Kamgarpour,
Ian Marquette, Henrik Schlichtkrull, and Ole Warnaar for helpful discussions and comments.

\appendix

\section{Layer sums and polynomials}
\label{Sec:LayerEx}

Here, we provide details of the layer sums and polynomials of the simple Lie algebras of rank $r\leq4$, as
well as $A_5$.
Because of well-known isomorphisms between the lower-rank Lie algebras, our focus will be on
\be
 A_1,\,A_2,\,A_3,\,A_4,\,A_5\qquad B_2,\,B_3,\,B_4,\qquad C_3,\,C_4,\qquad D_4,\qquad F_4,\qquad G_2.
\ee
The expressions for $A_1$ are trivially given by
\be
 \Lc_\la^{A_1}=\ch_\la,\qquad 
 R_{A_1}(\la)=R_{A_1}(\la_1)=1+\la_1,
\ee
while the expressions for $A_2$ and $G_2$ are given in (\ref{LA2}), (\ref{PA2}) and (\ref{LG2}), (\ref{PG2}),
respectively. The remaining examples are discussed in the following.

\subsection{The case $B_2$}
\label{Sec:Layer2}

For $B_2$, the number of non-simple positive roots is $|\Phi_+'|=2$. With the labelling convention
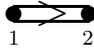
\setlength{\unitlength}{1mm}
\begin{center}
\begin{picture}(18,5)(-4.8,-2)
\thicklines
\put(0,0){\circle*{2}}
\put(10,0){\circle*{2}}
\put(0,-1){\line(1,0){10}}
\put(0,1){\line(1,0){10}}
\put(2.8,-1.5){{\LARGE $>$}}
\put(-0.7,-3.5){$_1$}
\put(9.1,-3.5){$_2$}
\end{picture}
\end{center}
the reduced Weyl vector is given by
\be
 \rho'=\tfrac{1}{2}(2\al_1+3\al_2),
\ee
while the layer sums and polynomial are given by
\be
 \Lc_\la^{B_2}=(\ch_\la+\ch_{\la-2\al_1-3\al_2})
 -(\ch_{\la-\al_1-\al_2}+\ch_{\la-\al_1-2\al_2})
\ee
and
\be
 R_{B_2}(\la)=1
  +2(\la_1+\la_2)
  +(2\la_1^2+\la_2^2+4\la_1\la_2).
\label{PB2}
\ee

\subsection{Rank-$3$ cases}
\label{Sec:Layer3}

For ease of comparison of the layer polynomials for $A_3$, $B_3$, and $C_3$, 
the $20$ distinct terms in $R_{X_3}$ are listed in the same order in the three cases.
Indeed, although simplifications are possible, no attempt has been made to
take into account the symmetries of the Dynkin diagrams.

For $A_3$, the number of non-simple positive roots is $|\Phi_+'|=3$. With the labelling convention
\setlength{\unitlength}{1mm}
\begin{center}
\begin{picture}(28,5)(-4.8,-2)
\thicklines
\put(0,0){\circle*{2}}
\put(10,0){\circle*{2}}
\put(20,0){\circle*{2}}
\put(0,0){\line(1,0){10}}
\put(10,0){\line(1,0){10}}
\put(-0.7,-3.5){$_1$}
\put(9.15,-3.5){$_2$}
\put(19.2,-3.5){$_3$}
\end{picture}
\end{center}
the reduced Weyl vector is given by
\be
 \rho'=\tfrac{1}{2}(2\al_1+3\al_2+2\al_3),
\ee
the layer sums by
\begin{align}
 \Lc_\la^{A_3}&=
 (\ch_\la-\ch_{\la-2\al_1-3\al_2-2\al_3})
 -(\ch_{\la-\al_1-\al_2}-\ch_{\la-\al_1-2\al_2-2\al_3})
 \nn
 &-(\ch_{\la-\al_2-\al_3}-\ch_{\la-2\al_1-2\al_2-\al_3})
 -(\ch_{\la-\al_1-\al_2-\al_3}-\ch_{\la-\al_1-2\al_2-\al_3}),
\end{align}
and the corresponding layer polynomial by
\begin{align}
 R_{A_3}(\la)&=1
  +\frac{1}{6}(11\la_1+14\la_2+11\la_3)
  +(\la_1^2+2\la_2^2+\la_3^2+4\la_1\la_2+4\la_2\la_3+3\la_1\la_3)
  +\frac{1}{6}(\la_1^3+4\la_2^3
  \nn
  &+\la_3^3+6\la_1^2\la_2+12\la_1\la_2^2+12\la_2^2\la_3+6\la_2\la_3^2
   +9\la_1^2\la_3+9\la_1\la_3^2+36\la_1\la_2\la_3).
\label{PA3}
\end{align}

For $B_3$, the number of non-simple positive roots is $|\Phi_+'|=6$.
With the labelling convention
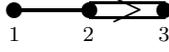
\setlength{\unitlength}{1mm}
\begin{center}
\begin{picture}(28,5)(-4.8,-2)
\thicklines
\put(0,0){\circle*{2}}
\put(10,0){\circle*{2}}
\put(20,0){\circle*{2}}
\put(0,0){\line(1,0){10}}
\put(10,-1){\line(1,0){10}}
\put(10,1){\line(1,0){10}}
\put(12.7,-1.5){{\LARGE $>$}}
\put(-0.7,-3.5){$_1$}
\put(9.15,-3.5){$_2$}
\put(19.2,-3.5){$_3$}
\end{picture}
\end{center}
the reduced Weyl vector is given by
\be
 \rho'=\tfrac{1}{2}(4\al_1+7\al_2+8\al_3),
\ee
the layer sums by
\begin{align}
 \Lc_\la^{B_3}&=
  (\ch_{\la}+\ch_{\la-4\al_1-7\al_2-8\al_3})
  -(\ch_{\la-\al_2-\al_3}+\ch_{\la-4\al_1-6\al_2-7\al_3})
  -(\ch_{\la-\al_1-\al_2}+\ch_{\la-3\al_1-6\al_2-8\al_3})
 \nn
 &
  -(\ch_{\la-\al_2-2\al_3}+\ch_{\la-4\al_1-6\al_2-6\al_3})
  +(\ch_{\la-2\al_2-3\al_3}+\ch_{\la-4\al_1-5\al_2-5\al_3})
 \nn
 &
  -(\ch_{\la-\al_1-\al_2-\al_3}+\ch_{\la-3\al_1-6\al_2-7\al_3})
  -(\ch_{\la-\al_1-\al_2-2\al_3}+\ch_{\la-3\al_1-6\al_2-6\al_3})
 \nn
 &
  +(\ch_{\la-\al_1-2\al_2-\al_3}+\ch_{\la-3\al_1-5\al_2-7\al_3})
  +(\ch_{\la-\al_1-2\al_2-2\al_3}+\ch_{\la-3\al_1-5\al_2-6\al_3})
 \nn
 &
  +2(\ch_{\la-\al_1-2\al_2-3\al_3}+\ch_{\la-3\al_1-5\al_2-5\al_3})
  +(\ch_{\la-\al_1-2\al_2-4\al_3}+\ch_{\la-3\al_1-5\al_2-4\al_3})
 \nn
 &
  -(\ch_{\la-\al_1-3\al_2-5\al_3}+\ch_{\la-3\al_1-4\al_2-3\al_3})
  -(\ch_{\la-\al_1-4\al_2-5\al_3}+\ch_{\la-3\al_1-3\al_2-3\al_3})
 \nn
 &
  +(\ch_{\la-2\al_1-2\al_2-\al_3}+\ch_{\la-2\al_1-5\al_2-7\al_3})
  +(\ch_{\la-2\al_1-2\al_2-2\al_3}+\ch_{\la-2\al_1-5\al_2-6\al_3})
 \nn
 &
  +(\ch_{\la-2\al_1-2\al_2-3\al_3}+\ch_{\la-2\al_1-5\al_2-5\al_3})
  -(\ch_{\la-2\al_1-3\al_2-3\al_3}+\ch_{\la-2\al_1-4\al_2-5\al_3})
 \nn
 &
  -(\ch_{\la-2\al_1-3\al_2-4\al_3}+\ch_{\la-2\al_1-4\al_2-4\al_3})
  -(\ch_{\la-2\al_1-3\al_2-5\al_3}+\ch_{\la-2\al_1-4\al_2-3\al_3}),
\end{align}
and the corresponding layer polynomial by
\begin{align}
 R_{B_3}(\la)&=1
  +\frac{1}{3}(8\la_1+10\la_2+9\la_3)
  +(2\la_1^2+8\la_2^2+3\la_3^2+8\la_1\la_2+12\la_2\la_3+6\la_1\la_3)
  +\frac{1}{3}(4\la_1^3+20\la_2^3
  \nn
  &+3\la_3^3+24\la_1^2\la_2+48\la_1\la_2^2+36\la_2^2\la_3+18\la_2\la_3^2
   +18\la_1^2\la_3+18\la_1\la_3^2+72\la_1\la_2\la_3).
\end{align}

For $C_3$, the number of non-simple positive roots is $|\Phi_+'|=6$. With the labelling convention
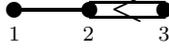
\setlength{\unitlength}{1mm}
\begin{center}
\begin{picture}(28,5)(-4.8,-2)
\thicklines
\put(0,0){\circle*{2}}
\put(10,0){\circle*{2}}
\put(20,0){\circle*{2}}
\put(0,0){\line(1,0){10}}
\put(10,-1){\line(1,0){10}}
\put(10,1){\line(1,0){10}}
\put(12.7,-1.5){{\LARGE $<$}}
\put(-0.7,-3.5){$_1$}
\put(9.15,-3.5){$_2$}
\put(19.2,-3.5){$_3$}
\end{picture}
\end{center}
the reduced Weyl vector is given by
\be
 \rho'=\tfrac{1}{2}(5\al_1+9\al_2+5\al_3),
\ee
the layer sums by
\begin{align}
 \Lc_\la^{C_3}&=
  (\ch_{\la}+\ch_{\la-5\al_1-9\al_2-5\al_3})
  -(\ch_{\la-\al_2-\al_3}+\ch_{\la-5\al_1-8\al_2-4\al_3})
  -(\ch_{\la-\al_1-\al_2}+\ch_{\la-4\al_1-8\al_2-5\al_3})
 \nn
 &
  -(\ch_{\la-2\al_2-\al_3}+\ch_{\la-5\al_1-7\al_2-4\al_3})
  -(\ch_{\la-\al_1-\al_2-\al_3}+\ch_{\la-4\al_1-8\al_2-4\al_3})
 \nn
 &
  +(\ch_{\la-3\al_2-2\al_3}+\ch_{\la-5\al_1-6\al_2-3\al_3})
  +(\ch_{\la-\al_1-2\al_2-2\al_3}+\ch_{\la-4\al_1-7\al_2-3\al_3})
 \nn
 &
  +(\ch_{\la-\al_1-3\al_2-\al_3}+\ch_{\la-4\al_1-6\al_2-4\al_3})
  +2(\ch_{\la-\al_1-3\al_2-2\al_3}+\ch_{\la-4\al_1-6\al_2-3\al_3})
 \nn
 &
  +(\ch_{\la-2\al_1-3\al_2-\al_3}+\ch_{\la-3\al_1-6\al_2-4\al_3})
  +(\ch_{\la-2\al_1-3\al_2-2\al_3}+\ch_{\la-3\al_1-6\al_2-3\al_3})
 \nn
 &
  +(\ch_{\la-3\al_1-3\al_2-\al_3}+\ch_{\la-2\al_1-6\al_2-4\al_3})
  -(\ch_{\la-\al_1-4\al_2-3\al_3}+\ch_{\la-4\al_1-5\al_2-2\al_3})
 \nn
 &
  -(\ch_{\la-2\al_1-4\al_2-2\al_3}+\ch_{\la-3\al_1-5\al_2-3\al_3})
  +(\ch_{\la-3\al_1-3\al_2-2\al_3}+\ch_{\la-2\al_1-6\al_2-3\al_3})
 \nn
 &
  -(\ch_{\la-\al_1-5\al_2-3\al_3}+\ch_{\la-4\al_1-4\al_2-2\al_3})
  -(\ch_{\la-2\al_1-4\al_2-3\al_3}+\ch_{\la-3\al_1-5\al_2-2\al_3})
 \nn
 &
  -(\ch_{\la-2\al_1-5\al_2-2\al_3}+\ch_{\la-3\al_1-4\al_2-3\al_3})
  -(\ch_{\la-3\al_1-4\al_2-2\al_3}+\ch_{\la-2\al_1-5\al_2-3\al_3}),
\end{align}
and the corresponding layer polynomial by
\begin{align}
 R_{C_3}(\la)&=1
  +\frac{1}{3}(7\la_1+11\la_2+9\la_3)
  +(2\la_1^2+5\la_2^2+6\la_3^2+8\la_1\la_2+12\la_2\la_3+6\la_1\la_3)
  +\frac{2}{3}(\la_1^3+5\la_2^3
  \nn
  &+6\la_3^3+6\la_1^2\la_2+12\la_1\la_2^2+18\la_2^2\la_3+18\la_2\la_3^2
   +9\la_1^2\la_3+18\la_1\la_3^2+36\la_1\la_2\la_3).
\end{align}

\subsection{Rank-$4$ cases}
\label{Sec:Layer4}

As the layer sums are rather involved for $r=4$, for $B_4$, $C_4$, $D_4$, and $F_4$, we only list the 
layer polynomials. For ease of comparison of the polynomials, 
the $70$ distinct terms in $R_{X_4}$ are listed in the same order in the five cases (including $A_4$).
Indeed, although simplifications are possible, no attempt has been made to
take into account the symmetries of the Dynkin diagrams.

For $A_4$, the number of non-simple positive roots is $|\Phi_+'|=6$. With the labelling convention
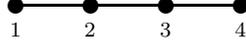
\setlength{\unitlength}{1mm}
\begin{center}
\begin{picture}(38,5)(-4.8,-2)
\thicklines
\put(0,0){\circle*{2}}
\put(10,0){\circle*{2}}
\put(20,0){\circle*{2}}
\put(30,0){\circle*{2}}
\put(0,0){\line(1,0){10}}
\put(10,0){\line(1,0){10}}
\put(20,0){\line(1,0){10}}
\put(-0.7,-3.5){$_1$}
\put(9.15,-3.5){$_2$}
\put(19.2,-3.5){$_3$}
\put(29.2,-3.5){$_4$}
\end{picture}
\end{center}
the reduced Weyl vector is given by
\be
 \rho'=\tfrac{1}{2}(3\al_1+5\al_2+5\al_3+3\al_4),
\ee
the layer sums by
\begin{align}
 \Lc_\la^{A_4}&=
  (\ch_\la+\ch_{\la-3\al_1-5\al_2-5\al_3-3\al_4})
  -(\ch_{\la-\al_1-\al_2}+\ch_{\la-2\al_1-4\al_2-5\al_3-3\al_4})
 \nn
 &
 -(\ch_{\la-\al_2-\al_3}+\ch_{\la-3\al_1-4\al_2-4\al_3-3\al_4})
 -(\ch_{\la-\al_3-\al_4}+\ch_{\la-3\al_1-5\al_2-4\al_3-2\al_4})
 \nn
 &
 -(\ch_{\la-\al_1-\al_2-\al_3}+\ch_{\la-2\al_1-4\al_2-4\al_3-3\al_4})
 -(\ch_{\la-\al_2-\al_3-\al_4}+\ch_{\la-3\al_1-4\al_2-4\al_3-2\al_4})
 \nn
 &
 +(\ch_{\la-\al_1-2\al_2-\al_3}+\ch_{\la-2\al_1-3\al_2-4\al_3-3\al_4})
 +(\ch_{\la-\al_2-2\al_3-\al_4}+\ch_{\la-3\al_1-4\al_2-3\al_3-2\al_4})
 \nn
 &
 +(\ch_{\la-2\al_1-2\al_2-\al_3}+\ch_{\la-\al_1-3\al_2-4\al_3-3\al_4})
 +(\ch_{\la-\al_2-2\al_3-2\al_4}+\ch_{\la-3\al_1-4\al_2-3\al_3-\al_4})
 \nn
 &
 +(\ch_{\la-\al_1-2\al_2-2\al_3}+\ch_{\la-2\al_1-3\al_2-3\al_3-3\al_4})
 +(\ch_{\la-2\al_2-2\al_3-\al_4}+\ch_{\la-3\al_1-3\al_2-3\al_3-2\al_4})
 \nn
 &
 +(\ch_{\la-\al_1-2\al_2-\al_3-\al_4}+\ch_{\la-2\al_1-3\al_2-4\al_3-2\al_4})
 +(\ch_{\la-\al_1-\al_2-2\al_3-\al_4}+\ch_{\la-2\al_1-4\al_2-3\al_3-2\al_4})
 \nn
 &
 +(\ch_{\la-2\al_1-2\al_2-\al_3-\al_4}+\ch_{\la-\al_1-3\al_2-4\al_3-2\al_4})
 +(\ch_{\la-\al_1-\al_2-2\al_3-2\al_4}+\ch_{\la-2\al_1-4\al_2-3\al_3-\al_4})
 \nn
 &
 +(\ch_{\la-\al_1-2\al_2-2\al_3-\al_4}+\ch_{\la-2\al_1-3\al_2-3\al_3-2\al_4})
 -(\ch_{\la-2\al_1-3\al_2-2\al_3}+\ch_{\la-\al_1-2\al_2-3\al_3-3\al_4})
 \nn
 &
 -(\ch_{\la-2\al_2-3\al_3-2\al_4}+\ch_{\la-3\al_1-3\al_2-2\al_3-\al_4})
 -(\ch_{\la-\al_1-3\al_2-2\al_3-\al_4}+\ch_{\la-2\al_1-2\al_2-3\al_3-2\al_4})
 \nn
 &
 -(\ch_{\la-\al_1-2\al_2-3\al_3-\al_4}+\ch_{\la-2\al_1-3\al_2-2\al_3-2\al_4})
 -2(\ch_{\la-2\al_1-3\al_2-2\al_3-\al_4}+\ch_{\la-\al_1-2\al_2-3\al_3-2\al_4})
 \nn
 &
 -(\ch_{\la-\al_1-3\al_2-3\al_3-\al_4}+\ch_{\la-2\al_1-2\al_2-2\al_3-2\al_4}),
\end{align}
and the corresponding layer polynomial by
\begin{align}
 R_{A_4}(\la)&=1
  +\frac{5}{12}(5\la_1+7\la_2+7\la_3+5\la_4)
  +\frac{5}{24}(7\la_1^2+17\la_2^2+17\la_3^2+7\la_4^2+28\la_1\la_2+40\la_2\la_3
  \nn
  &+28\la_3\la_4+26\la_1\la_3+26\la_2\la_4+20\la_1\la_4)
  +\frac{5}{12}(\la_1^3+5\la_2^3+5\la_3^3+\la_4^3+6\la_1^2\la_2+12\la_1\la_2^2
  \nn
  &+18\la_2^2\la_3+18\la_2\la_3^2+12\la_3^2\la_4+6\la_3\la_4^2+9\la_1^2\la_3+15\la_1\la_3^2+15\la_2^2\la_4
   +9\la_2\la_4^2+6\la_1^2\la_4+6\la_1\la_4^2
  \nn
  &+36\la_1\la_2\la_3+36\la_2\la_3\la_4+24\la_1\la_2\la_4+24\la_1\la_3\la_4)
  +\frac{1}{24}(\la_1^4+11\la_2^4+11\la_3^4+\la_4^4+8\la_1^3\la_2
  \nn
  &+32\la_1\la_2^3+56\la_2^3\la_3+56\la_2\la_3^3+32\la_3^3\la_4+8\la_3\la_4^3
   +12\la_1^3\la_3+68\la_1\la_3^3+68\la_2^3\la_4+12\la_2\la_4^3
  \nn
  &+16\la_1^3\la_4+16\la_1\la_4^3+24\la_1^2\la_2^2+96\la_2^2\la_3^2
    +24\la_3^2\la_4^2+54\la_1^2\la_3^2+54\la_2^2\la_4^2+36\la_1^2\la_4^2
  \nn
  &+72\la_1^2\la_2\la_3+144\la_1\la_2^2\la_3+216\la_1\la_2\la_3^2+216\la_2^2\la_3\la_4+144\la_2\la_3^2\la_4
   +72\la_2\la_3\la_4^2+96\la_1^2\la_2\la_4
  \nn
  &+192\la_1\la_2^2\la_4+144\la_1\la_2\la_4^2
   +144\la_1^2\la_3\la_4+192\la_1\la_3^2\la_4+96\la_1\la_3\la_4^2+576\la_1\la_2\la_3\la_4).
\end{align}

For $B_4$, the number of non-simple positive roots is $|\Phi_+'|=12$. With the labelling convention
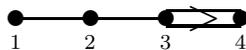
\setlength{\unitlength}{1mm}
\begin{center}
\begin{picture}(38,5)(-4.8,-2)
\thicklines
\put(0,0){\circle*{2}}
\put(10,0){\circle*{2}}
\put(20,0){\circle*{2}}
\put(30,0){\circle*{2}}
\put(0,0){\line(1,0){10}}
\put(10,0){\line(1,0){10}}
\put(20,-1){\line(1,0){10}}
\put(20,1){\line(1,0){10}}
\put(22.7,-1.5){{\LARGE $>$}}
\put(-0.7,-3.5){$_1$}
\put(9.15,-3.5){$_2$}
\put(19.2,-3.5){$_3$}
\put(29.2,-3.5){$_4$}
\end{picture}
\end{center}
the reduced Weyl vector is given by
\be
 \rho'=\tfrac{1}{2}(6\al_1+11\al_2+14\al_3+15\al_4),
\ee
while the layer polynomial is given by
\begin{align}
 R_{B_4}(\la)&=1
  +\frac{4}{3}(2\la_1+4\la_2+3\la_3+3\la_4)
  +\frac{2}{3}(5\la_1^2+12\la_2^2+25\la_3^2+9\la_4^2+20\la_1\la_2+28\la_2\la_3+36\la_3\la_4
  \nn
  &+18\la_1\la_3+20\la_2\la_4+16\la_1\la_4)
  +\frac{4}{3}(\la_1^3+8\la_2^3+21\la_3^3+3\la_4^3+6\la_1^2\la_2+12\la_1\la_2^2+36\la_2^2\la_3
  \nn
  &+54\la_2\la_3^2+36\la_3^2\la_4+18\la_3\la_4^2+9\la_1^2\la_3+27\la_1\la_3^2+24\la_2^2\la_4+18\la_2\la_4^2
    +6\la_1^2\la_4+9\la_1\la_4^2
  \nn
  &+36\la_1\la_2\la_3+72\la_2\la_3\la_4+24\la_1\la_2\la_4+36\la_1\la_3\la_4)
  +\frac{1}{3}(2\la_1^4+24\la_2^4+46\la_3^4+3\la_4^4+16\la_1^3\la_2
  \nn
  &+64\la_1\la_2^3+128\la_2^3\la_3+176\la_2\la_3^3+96\la_3^3\la_4+24\la_3\la_4^3
   +24\la_1^3\la_3+168\la_1\la_3^3+80\la_2^3\la_4+24\la_2\la_4^3
  \nn
  &+16\la_1^3\la_4+24\la_1\la_4^3+48\la_1^2\la_2^2+240\la_2^2\la_3^2
    +72\la_3^2\la_4^2+108\la_1^2\la_3^2+72\la_2^2\la_4^2+36\la_1^2\la_4^2
  \nn
  &+144\la_1^2\la_2\la_3+288\la_1\la_2^2\la_3+432\la_1\la_2\la_3^2+288\la_2^2\la_3\la_4
   +288\la_2\la_3^2\la_4+144\la_2\la_3\la_4^2+96\la_1^2\la_2\la_4
  \nn
  &+192\la_1\la_2^2\la_4+144\la_1\la_2\la_4^2
   +144\la_1^2\la_3\la_4+288\la_1\la_3^2\la_4+144\la_1\la_3\la_4^2+576\la_1\la_2\la_3\la_4).
\end{align}

For $C_4$, the number of non-simple positive roots is $|\Phi_+'|=12$. With the labelling convention
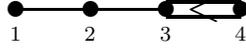
\setlength{\unitlength}{1mm}
\begin{center}
\begin{picture}(38,5)(-4.8,-2)
\thicklines
\put(0,0){\circle*{2}}
\put(10,0){\circle*{2}}
\put(20,0){\circle*{2}}
\put(30,0){\circle*{2}}
\put(0,0){\line(1,0){10}}
\put(10,0){\line(1,0){10}}
\put(20,-1){\line(1,0){10}}
\put(20,1){\line(1,0){10}}
\put(22.7,-1.5){{\LARGE $<$}}
\put(-0.7,-3.5){$_1$}
\put(9.15,-3.5){$_2$}
\put(19.2,-3.5){$_3$}
\put(29.2,-3.5){$_4$}
\end{picture}
\end{center}
the reduced Weyl vector is given by
\be
 \rho'=\tfrac{1}{2}(7\al_1+13\al_2+17\al_3+9\al_4),
\ee
while the layer polynomial is given by
\begin{align}
 R_{C_4}(\la)&=1
  +\frac{4}{3}(2\la_1+3\la_2+4\la_3+3\la_4)
  +\frac{2}{3}(4\la_1^2+12\la_2^2+17\la_3^2+18\la_4^2+16\la_1\la_2+32\la_2\la_3+36\la_3\la_4
  \nn
  &+18\la_1\la_3+20\la_2\la_4+16\la_1\la_4)
  +\frac{4}{3}(\la_1^3+6\la_2^3+11\la_3^3+12\la_4^3+6\la_1^2\la_2+12\la_1\la_2^2+24\la_2^2\la_3
  \nn
  &+30\la_2\la_3^2+36\la_3^2\la_4+36\la_3\la_4^2+9\la_1^2\la_3+18\la_1\la_3^2+24\la_2^2\la_4+36\la_2\la_4^2
    +6\la_1^2\la_4+18\la_1\la_4^2
  \nn
  &+36\la_1\la_2\la_3+72\la_2\la_3\la_4+24\la_1\la_2\la_4+36\la_1\la_3\la_4)
  +\frac{1}{3}(\la_1^4+12\la_2^4+23\la_3^4+24\la_4^4+8\la_1^3\la_2
  \nn
  &+32\la_1\la_2^3+64\la_2^3\la_3+88\la_2\la_3^3+96\la_3^3\la_4+96\la_3\la_4^3
   +12\la_1^3\la_3+84\la_1\la_3^3+80\la_2^3\la_4+96\la_2\la_4^3
  \nn
  &+16\la_1^3\la_4+96\la_1\la_4^3+24\la_1^2\la_2^2+120\la_2^2\la_3^2
    +144\la_3^2\la_4^2+54\la_1^2\la_3^2+144\la_2^2\la_4^2+72\la_1^2\la_4^2
  \nn
  &+72\la_1^2\la_2\la_3+144\la_1\la_2^2\la_3+216\la_1\la_2\la_3^2+288\la_2^2\la_3\la_4+288\la_2\la_3^2\la_4
   +288\la_2\la_3\la_4^2+96\la_1^2\la_2\la_4
  \nn
  &+192\la_1\la_2^2\la_4+288\la_1\la_2\la_4^2
   +144\la_1^2\la_3\la_4+288\la_1\la_3^2\la_4+288\la_1\la_3\la_4^2+576\la_1\la_2\la_3\la_4).
\end{align}

For $D_4$, the number of non-simple positive roots is $|\Phi_+'|=8$. With the labelling convention
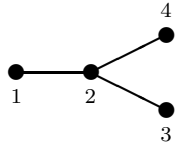
\setlength{\unitlength}{1mm}
\begin{center}
\begin{picture}(20,15)(0,-7)
\thicklines
\put(0,0){\circle*{2}}
\put(10,0){\circle*{2}}
\put(20,-5){\circle*{2}}
\put(20,5){\circle*{2}}
\put(0,0){\line(1,0){10}}
\put(10,0){\line(2,1){10}}
\put(10,0){\line(2,-1){10}}
\put(-0.7,-3.5){$_1$}
\put(9.15,-3.5){$_2$}
\put(19.2,-8.5){$_3$}
\put(19.2,8){$_4$}
\end{picture}
\end{center}
the reduced Weyl vector is given by
\be
 \rho'=\tfrac{1}{2}(5\al_1+9\al_2+5\al_3+5\al_4),
\ee
while the layer polynomial is given by
\begin{align}
 R_{D_4}(\la)&=1
  +\frac{4}{3}(2\la_1+3\la_2+2\la_3+2\la_4)
  +\frac{2}{3}(4\la_1^2+12\la_2^2+4\la_3^2+4\la_4^2+16\la_1\la_2+16\la_2\la_3
  \nn
  &+9\la_3\la_4+9\la_1\la_3+16\la_2\la_4+9\la_1\la_4)
  +\frac{2}{3}(2\la_1^3+12\la_2^3+2\la_3^3+2\la_4^3+12\la_1^2\la_2+24\la_1\la_2^2
  \nn
  &+24\la_2^2\la_3+12\la_2\la_3^2+9\la_3^2\la_4+9\la_3\la_4^2+9\la_1^2\la_3+9\la_1\la_3^2+24\la_2^2\la_4
   +12\la_2\la_4^2+9\la_1^2\la_4+9\la_1\la_4^2
  \nn
  &+36\la_1\la_2\la_3+36\la_2\la_3\la_4+36\la_1\la_2\la_4+18\la_1\la_3\la_4)
  +\frac{1}{3}(\la_1^4+12\la_2^4+\la_3^4+\la_4^4+8\la_1^3\la_2
  \nn
  &+32\la_1\la_2^3+32\la_2^3\la_3+8\la_2\la_3^3+6\la_3^3\la_4+6\la_3\la_4^3
   +6\la_1^3\la_3+6\la_1\la_3^3+32\la_2^3\la_4+8\la_2\la_4^3
  \nn
  &+6\la_1^3\la_4+6\la_1\la_4^3+24\la_1^2\la_2^2+24\la_2^2\la_3^2
    +9\la_3^2\la_4^2+9\la_1^2\la_3^2+24\la_2^2\la_4^2+9\la_1^2\la_4^2
  \nn
  &+36\la_1^2\la_2\la_3+72\la_1\la_2^2\la_3+36\la_1\la_2\la_3^2+72\la_2^2\la_3\la_4+36\la_2\la_3^2\la_4
   +36\la_2\la_3\la_4^2+36\la_1^2\la_2\la_4
  \nn
  &+72\la_1\la_2^2\la_4+36\la_1\la_2\la_4^2
   +36\la_1^2\la_3\la_4+36\la_1\la_3^2\la_4+36\la_1\la_3\la_4^2+144\la_1\la_2\la_3\la_4).
\end{align}

For $F_4$, the number of non-simple positive roots is $|\Phi_+'|=20$. With the labelling convention
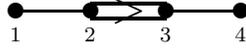
\setlength{\unitlength}{1mm}
\begin{center}
\begin{picture}(38,5)(-4.8,-2)
\thicklines
\put(0,0){\circle*{2}}
\put(10,0){\circle*{2}}
\put(20,0){\circle*{2}}
\put(30,0){\circle*{2}}
\put(0,0){\line(1,0){10}}
\put(10,-1){\line(1,0){10}}
\put(10,1){\line(1,0){10}}
\put(20,0){\line(1,0){10}}
\put(12.7,-1.5){{\LARGE $>$}}
\put(-0.7,-3.5){$_1$}
\put(9.15,-3.5){$_2$}
\put(19.2,-3.5){$_3$}
\put(29.2,-3.5){$_4$}
\end{picture}
\end{center}
the reduced Weyl vector is given by
\be
 \rho'=\tfrac{1}{2}(15\al_1+29\al_2+41\al_3+21\al_4),
\ee
while the layer polynomial is given by
\begin{align}
 R_{F_4}(\la)&=1
  +4(2\la_1+\la_2+2\la_3+\la_4)
  +2(4\la_1^2+22\la_2^2+13\la_3^2+4\la_4^2+16\la_1\la_2+36\la_2\la_3+16\la_3\la_4
  \nn
  &+16\la_1\la_3+8\la_2\la_4+12\la_1\la_4)
  +4(4\la_1^3+32\la_2^3+13\la_3^3+2\la_4^3+24\la_1^2\la_2+48\la_1\la_2^2+72\la_2^2\la_3
  \nn
  &+54\la_2\la_3^2+24\la_3^2\la_4+12\la_3\la_4^2+18\la_1^2\la_3+24\la_1\la_3^2+42\la_2^2\la_4
   +18\la_2\la_4^2+6\la_1^2\la_4+6\la_1\la_4^2
  \nn
  &+72\la_1\la_2\la_3+72\la_2\la_3\la_4+24\la_1\la_2\la_4+24\la_1\la_3\la_4)
  +2(8\la_1^4+116\la_2^4+29\la_3^4+2\la_4^4+64\la_1^3\la_2
  \nn
  &+256\la_1\la_2^3+336\la_2^3\la_3+168\la_2\la_3^3+64\la_3^3\la_4+16\la_3\la_4^3
   +48\la_1^3\la_3+104\la_1\la_3^3+208\la_2^3\la_4+24\la_2\la_4^3
  \nn
  &+32\la_1^3\la_4+16\la_1\la_4^3+192\la_1^2\la_2^2+360\la_2^2\la_3^2
    +48\la_3^2\la_4^2+108\la_1^2\la_3^2+108\la_2^2\la_4^2+36\la_1^2\la_4^2
  \nn
  &+288\la_1^2\la_2\la_3+576\la_1\la_2^2\la_3+432\la_1\la_2\la_3^2+432\la_2^2\la_3\la_4
   +288\la_2\la_3^2\la_4+144\la_2\la_3\la_4^2+192\la_1^2\la_2\la_4
  \nn
  &+384\la_1\la_2^2\la_4+144\la_1\la_2\la_4^2
   +144\la_1^2\la_3\la_4+192\la_1\la_3^2\la_4+96\la_1\la_3\la_4^2+576\la_1\la_2\la_3\la_4).
\end{align}

\subsection{The case $A_5$}
\label{Sec:LayerA5}

For $A_5$, the number of non-simple positive roots is $|\Phi_+'|=10$. With the labelling convention
\setlength{\unitlength}{1mm}
\begin{center}
\begin{picture}(48,5)(-4.8,-2)
\thicklines
\put(0,0){\circle*{2}}
\put(10,0){\circle*{2}}
\put(20,0){\circle*{2}}
\put(30,0){\circle*{2}}
\put(40,0){\circle*{2}}
\put(0,0){\line(1,0){10}}
\put(10,0){\line(1,0){10}}
\put(20,0){\line(1,0){10}}
\put(30,0){\line(1,0){10}}
\put(-0.7,-3.5){$_1$}
\put(9.15,-3.5){$_2$}
\put(19.2,-3.5){$_3$}
\put(29.2,-3.5){$_4$}
\put(39.2,-3.5){$_5$}
\end{picture}
\end{center}
the reduced Weyl vector is given by
\be
 \rho'=\tfrac{1}{2}(4\al_1+7\al_2+8\al_3+7\al_4+4\al_5),
\ee
while the layer polynomial is given by
\begin{align}
 R_{A_5}(\la)&=1
 +\frac{1}{60}(137\la_1+202\la_2+222\la_3+202\la_4+137\la_5)
 +\frac{5}{8}(3\la_1^2+8\la_2^2+10\la_3^2+8\la_4^2+3\la_5^2
 \nn 
 &+12\la_1\la_2+20\la_2\la_3+20\la_3\la_4+12\la_4\la_5+12\la_1\la_3+16\la_2\la_4+12\la_3\la_5+12\la_1\la_4
  +12\la_2\la_5
 \nn
 &+8\la_1\la_5)
 +\frac{1}{24}(17\la_1^3+94\la_2^3+138\la_3^3+94\la_4^3+17\la_5^3+102\la_1^2\la_2+204\la_1\la_2^2
   +360\la_2^2\la_3
 \nn
 &+414\la_2\la_3^2+414\la_3^2\la_4+360\la_3\la_4^2+204\la_4^2\la_5+102\la_4\la_5^2
   +153\la_1^2\la_3+294\la_1\la_3^2+358\la_2^2\la_4
 \nn
 &+358\la^2\la_4^2+294\la_3^2\la_5+153\la_3\la_5^2+
   +124\la_1^2\la_4+236\la_1\la_4^2+236\la_2^2\la_5+124\la_2\la_5^2+95\la_1\la_5
 \nn
 &+95\la_1\la_5^2+612\la_1\la_2\la_3+936\la_2\la_3\la_4+612\la_3\la_4\la_5
   +496\la_1\la_2\la_4+564\la_1\la_3\la_4+564\la_2\la_3\la_5
 \nn
 &+496\la_2\la_4\la_5+380\la_1\la_2\la_5+380\la_1\la_4\la_5+360\la_1\la_3\la_5)
 +p_4+p_5,
\end{align}
where
\begin{align}
 p_4&=\frac{1}{8}(\la_1^4+12\la_2^4+22\la_3^4+12\la_4^4+\la_5^4)
  +\frac{1}{4}(4\la_1^3\la_2+16\la_1\la_2^3+32\la_2^3\la_3+44\la_2\la_3^3+44\la_3^3\la_4
 \nn
 &+32\la_3\la_4^3+16\la_4^3\la_5+4\la_4\la_5^3+6\la_1^3\la_3+39\la_1\la_3^3+40\la_2^3\la_4+40\la_2\la_4^3
  +39\la_3^3\la_5+6\la_3\la_5^3+8\la_1^3\la_4
 \nn
 &+28\la_1\la_4^3+28\la_2^3\la_5+8\la_2\la_5^3+5\la_1^3\la_5+5\la_1\la_5^3+12\la_1^2\la_2^2
  +60\la_2^2\la_3^2+60\la_3^2\la_4^2+12\la_4^2\la_5^2+27\la_1^2\la_3^2
 \nn
 &+64\la_2^2\la_4^2+27\la_3^2\la_5^2+28\la_1^2\la_4^2+28\la_2^2\la_5^2+10\la_1^2\la_5^2)
  +(9\la_1^2\la_2\la_3+18\la_1\la_2^2\la_3+27\la_1\la_2\la_3^2
 \nn
 &+36\la_2^2\la_3\la_4+36\la_2\la_3^2\la_4+36\la_2\la_3\la_4^2+27\la_3^2\la_4\la_5+18\la_3\la_4^2\la_5
  +9\la_3\la_4\la_5^2+12\la_1^2\la_2\la_4
 \nn
 &+24\la_1\la_2^2\la_4+18\la_1^2\la_3\la_4)
  +\frac{1}{2}(63\la_1\la_3^2\la_4+56\la_1\la_2\la_4^2+54\la_1\la_3\la_4^2+54\la_2^2\la_3\la_5
  +63\la_2\la_3^2\la_5
 \nn
 &+56\la_2^2\la_4\la_5+48\la_2\la_4^2\la_5+36\la_2\la_3\la_5^2+24\la_2\la_4\la_5^2)
  +\frac{1}{4}(30\la_1^2\la_2\la_5
  +60\la_1\la_2^2\la_5+40\la_1^2\la_4\la_5
 \nn
 &+60\la_1\la_4^2\la_5+40\la_1\la_2\la_5^2+30\la_1\la_4\la_5^2+45\la_1^2\la_3\la_5+90\la_1\la_3^2\la_5
  +45\la_1\la_3\la_5^2)+(72\la_1\la_2\la_3\la_4
 \nn
 &+72\la_2\la_3\la_4\la_5+45\la_1\la_2\la_3\la_5+45\la_1\la_3\la_4\la_5
   +40\la_1\la_2\la_4\la_5)
\end{align}
and
\begin{align}
 p_5&=\frac{1}{120}(\la_1^5+26\la_2^5+66\la_3^5+26\la_4^5+\la_5^5)
  +\frac{1}{24}(2\la_1^4\la_2+16\la_1\la_2^4+36\la_2^4\la_3+66\la_2\la_3^4+66\la_3^4\la_4
 \nn
 &+36\la_3\la_4^4+16\la_4^4\la_5+2\la_4\la_5^4+3\la_1^4\la_3+66\la_1\la_3^4+46\la_2^4\la_4+46\la_2\la_4^4
  +66\la_3^4\la_5+3\la_3\la_5^4+4\la_1^4\la_4
 \nn
 &+56\la_1\la_4^4+56\la_2^4\la_5+4\la_2\la_5^4+5\la_1^4\la_5+5\la_1\la_5^4+8\la_1^3\la_2^2
  +16\la_1^2\la_2^3+96\la_2^3\la_3^2+120\la_2^2\la_3^3+120\la_3^3\la_4^2
 \nn
 &+96\la_3^2\la_4^3+16\la_4^3\la_5^2+8\la_4^2\la_5^3+18\la_1^3\la_3^2+54\la_1^2\la_3^3+148\la_2^3\la_4^2
  +148\la_2^2\la_4^3+54\la_3^3\la_5^2+18\la_3^2\la_5^3
 \nn
 &+32\la_1^3\la_4^2+88\la_1^2\la_4^3+88\la_2^3\la_5^2+32\la_2^2\la_5^3+20\la_1^3\la_5^2
  +20\la_1^2\la_5^3)
 +\frac{1}{12}(12\la_1^3\la_2\la_3+48\la_1\la_2^3\la_3
 \nn
 &+108\la_1\la_2\la_3^3+120\la_2^3\la_3\la_4+144\la_2\la_3^3\la_4+120\la_2\la_3\la_4^3+108\la_3^3\la_4\la_5
  +48\la_3\la_4^3\la_5+12\la_3\la_4\la_5^3
 \nn
 &+16\la_1^3\la_2\la_4+64\la_1\la_2^3\la_4+24\la_1^3\la_3\la_4+156\la_1\la_3^3\la_4
  +144\la_1\la_3\la_4^3+176\la_1\la_2\la_4^3+144\la_2^3\la_3\la_5
 \nn
 &+156\la_2\la_3^3\la_5+176\la_2^3\la_4\la_5+64\la_2\la_4^3\la_5+24\la_2\la_3\la_5^3+16\la_2\la_4\la_5^3
  +20\la_1^3\la_2\la_5+80\la_1\la_2^3\la_5+40\la_1^3\la_4\la_5
 \nn
 &+80\la_1\la_4^3\la_5+40\la_1\la_2\la_5^3+20\la_1\la_4\la_5^3+30\la_1^3\la_3\la_5
  +180\la_1\la_3^3\la_5+30\la_1\la_3\la_5^3+36\la_1^2\la_2^2\la_3+54\la_1^2\la_2\la_3^2
 \nn
 &+108\la_1\la_2^2\la_3^2
  +216\la_2^2\la_3^2\la_4+252\la_2^2\la_3\la_4^2+216\la_2\la_3^2\la_4^2+108\la_3^2\la_4^2\la_5
  +54\la_3^2\la_4\la_5^2+36\la_3\la_4^2\la_5^2
 \nn
 &+48\la_1^2\la_2^2\la_4+108\la_1^2\la_3^2\la_4+96\la_1^2\la_2\la_4^2+192\la_1\la_2^2\la_4^2
  +144\la_1^2\la_3\la_4^2+252\la_1\la_3^2\la_4^2+192\la_2^2\la_4^2\la_5
 \nn
 &+252\la_2^2\la_3^2\la_5+144\la_2^2\la_3\la_5^2+108\la_2\la_3^2\la_5^2
  +96\la_2^2\la_4\la_5^2+48\la_2\la_4^2\la_5^2+60\la_1^2\la_2^2\la_5+120\la_1^2\la_4^2\la_5
 \nn
 &+60\la_1^2\la_2\la_5^2+120\la_1\la_2^2\la_5^2+60\la_1^2\la_4\la_5^2
  +60\la_1\la_4^2\la_5^2+135\la_1^2\la_3^2\la_5+90\la_1^2\la_3\la_5^2+135\la_1\la_3^2\la_5^2)
 \nn
 &+(12\la_1^2\la_2\la_3\la_4+24\la_1\la_2^2\la_3\la_4+36\la_1\la_2\la_3^2\la_4
  +48\la_1\la_2\la_3\la_4^2+48\la_2^2\la_3\la_4\la_5+36\la_2\la_3^2\la_4\la_5
 \nn
 &+24\la_2\la_3\la_4^2\la_5+12\la_2\la_3\la_4\la_5^2+15\la_1^2\la_2\la_3\la_5
  +30\la_1\la_2^2\la_3\la_5+45\la_1\la_2\la_3^2\la_5+20\la_1^2\la_2\la_4\la_5
 \nn 
 &+40\la_1\la_2^2\la_4\la_5+30\la_1^2\la_3\la_4\la_5+45\la_1\la_3^2\la_4\la_5
  +40\la_1\la_2\la_4^2\la_5+30\la_1\la_3\la_4^2\la_5+30\la_1\la_2\la_3\la_5^2
 \nn
 &+20\la_1\la_2\la_4\la_5^2+15\la_1\la_3\la_4\la_5^2+120\la_1\la_2\la_3\la_4\la_5).
\end{align}

\section{Discrete polytope volumes}
\label{Sec:Equiv}

Let $\la,\mu\in P_+$. The condition (\ref{lamuA}) for $\mu\in P(\la)$ means that
\be
 \la-\mu=n_1\al_1+\ldots+n_r\al_r
\label{lamu}
\ee
for some $n_1,\ldots,n_r\in\mathbb{N}_0$. Below, we use this to evaluate the sum of orbit lengths in (\ref{Ga}).
To simplify the characterisation of the various orbit lengths, 
we shall use a notation where $\nu_1,\nu_2,\nu_3\in\mathbb{N}$.

\subsection{Rank-$2$ cases}
\label{Sec:Equiv2}

For $A_2$,
\be
 |O_{\nu_1\om_1+\nu_2\om_2}|=6,\qquad |O_{\nu_1\om_1}|=|O_{\nu_2\om_2}|=3,\qquad |O_{0}|=1,
\ee 
and the condition (\ref{lamu}) requires
\be
 \la_1-2n_1+n_2,\,\la_2+n_1-2n_2\in\mathbb{N}_0.
\ee
This implies that
\be
 \sum_{\mu\in P_+(\la)}\!\!|O_\mu|
 =\sum_{n_1=0}^{\big\lfloor\frac{2\la_1+\la_2-3}{3}\big\rfloor}
  \sum_{n_2=\max(0,2n_1-\la_1+1)}^{\big\lfloor\frac{n_1+\la_2-1}{2}\big\rfloor} \!\!6
 +\sum_{n_1=\big\lceil\frac{\la_1}{2}\big\rceil}^{\big\lfloor\frac{2\la_1+\la_2-1}{3}\big\rfloor} \!\!3
 +\sum_{n_2=\big\lceil\frac{\la_2}{2}\big\rceil}^{\big\lfloor\frac{\la_1+2\la_2-1}{3}\big\rfloor} \!\!3
 +\sum_{n_1=\big\lceil\frac{2\la_1+\la_2}{3}\big\rceil}^{\big\lfloor\frac{2\la_1+\la_2}{3}\big\rfloor}
  \sum_{n_2=\big\lceil\frac{\la_1+2\la_2}{3}\big\rceil}^{\big\lfloor\frac{\la_1+2\la_2}{3}\big\rfloor} \!\!1,
\ee
which is seen to agree with (\ref{PA2}).

For $B_2$,
\be
 |O_{\nu_1\om_1+\nu_2\om_2}|=8,\qquad |O_{\nu_1\om_1}|=|O_{\nu_2\om_2}|=4,\qquad |O_{0}|=1,
\ee 
and the condition (\ref{lamu}) requires
\be
 \la_1-2n_1+n_2,\,\la_2+2n_1-2n_2\in\mathbb{N}_0.
\ee
This implies that
\be
 \sum_{\mu\in P_+(\la)}\!\!|O_\mu|
 =\sum_{n_1=0}^{\big\lfloor\frac{2\la_1+\la_2-3}{2}\big\rfloor}
  \sum_{n_2=\max(0,2n_1-\la_1+1)}^{\big\lfloor\frac{2n_1+\la_2-1}{2}\big\rfloor} \!\!8
 +\sum_{n_1=\big\lceil\frac{\la_1}{2}\big\rceil}^{\big\lfloor\frac{2\la_1+\la_2-1}{2}\big\rfloor} \!\!4
 +\sum_{n_2=\big\lceil\frac{\la_2}{2}\big\rceil}^{\la_1+\la_2-1}
  \sum_{n_1=\big\lceil\frac{2n_2-\la_2}{2}\big\rceil}^{\big\lfloor\frac{2n_2-\la_2}{2}\big\rfloor} \!\!4
 +\sum_{n_1=\big\lceil\frac{2\la_1+\la_2}{2}\big\rceil}^{\big\lfloor\frac{2\la_1+\la_2}{2}\big\rfloor} \!\!1,
\ee
which is seen to agree with (\ref{PB2}).

For $G_2$,
\be
 |O_{\nu_1\om_1+\nu_2\om_2}|=12,\qquad |O_{\nu_1\om_1}|=|O_{\nu_2\om_2}|=6,\qquad |O_{0}|=1,
\ee
and the condition (\ref{lamu}) requires
\be
 \la_1-2n_1+n_2,\,\la_2+3n_1-2n_2\in\mathbb{N}_0.
\ee
This implies that
\be
 \sum_{\mu\in P_+(\la)}\!\!|O_\mu|
 =\sum_{n_1=0}^{2\la_1+\la_2-3}
  \sum_{n_2=\max(0,2n_1-\la_1+1)}^{\big\lfloor\frac{3n_1+\la_2-1}{2}\big\rfloor} \!\!12
 +\sum_{n_1=\big\lceil\frac{\la_1}{2}\big\rceil}^{2\la_1+\la_2-1} \!\!6
 +\sum_{n_2=\big\lceil\frac{\la_2}{2}\big\rceil}^{3\la_1+2\la_2-3}
  \sum_{n_1=\big\lceil\frac{2n_2-\la_2}{3}\big\rceil}^{\big\lfloor\frac{2n_2-\la_2}{3}\big\rfloor} \!\!6
 +1,
\ee
which is seen to agree with (\ref{PG2}). 
This computation also explains the $R_{G_2}(\la)$ property observed immediately following (\ref{PG2mod6}).

\subsection{The case $A_3$}
\label{Sec:EquivA3}

For $A_3$,
\be
\begin{array}{c}
 |O_{\nu_1\om_1+\nu_2\om_2+\nu_3\om_3}|=24,\qquad 
 |O_{\nu_1\om_1+\nu_2\om_2}|=|O_{\nu_1\om_1+\nu_3\om_3}|=|O_{\nu_2\om_2+\nu_3\om_3}|=12,
 \\[.3cm]
 |O_{\nu_2\om_2}|=6,\qquad |O_{\nu_1\om_1}|=|O_{\nu_3\om_3}|=4,\qquad |O_{0}|=1,
\end{array}
\ee
and the condition (\ref{lamu}) requires
\be
 \la_1-2n_1+n_2,\,\la_2+n_1-2n_2+n_3,\,\la_3+n_2-2n_3\in\mathbb{N}_0.
\ee
This implies that
\begin{align}
 \sum_{\mu\in P_+(\la)}\!\!|O_\mu|
 &=\sum_{n_1=0}^{\big\lfloor\frac{3\la_1+2\la_2+\la_3-6}{4}\big\rfloor}
  \sum_{n_3=0}^{\big\lfloor\frac{\la_1+2\la_2+3\la_3-6}{4}\big\rfloor}
  \sum_{n_2=\max(0,1-\la_1+2n_1,1-\la_3+2n_3)}^{\big\lfloor\frac{\la_2+n_1+n_3-1}{2}\big\rfloor} \!\!24
 \nonumber\\[.2cm]
 &+\sum_{n_3=\big\lceil\frac{\la_3}{2}\big\rceil}^{\big\lfloor\frac{\la_1+2\la_2+3\la_3-3}{4}\big\rfloor}
  \sum_{n_1=\max(0,1-\la_2-2\la_3+3n_3)}^{\big\lfloor\frac{\la_1-\la_3+2n_3-1}{2}\big\rfloor} \!\!12
 +\sum_{n_1=\big\lceil\frac{\la_1}{2}\big\rceil}^{\big\lfloor\frac{3\la_1+2\la_2+\la_3-3}{4}\big\rfloor}
  \sum_{n_3=\max(0,1-2\la_1-\la_2+3n_1)}^{\big\lfloor\frac{-\la_1+\la_3+2n_1-1}{2}\big\rfloor} \!\!12
 \nonumber\\[.2cm]
 & +\sum_{n_2=0}^{\big\lfloor\frac{\la_1+2\la_2+\la_3-2}{2}\big\rfloor}
  \sum_{n_1=\max(0,\big\lceil\frac{1-2\la_2-\la_3+3n_2}{2}\big\rceil)}^{
    \min(-\la_2+2n_2,\big\lfloor\frac{\la_1+n_2-1}{2}\big\rfloor)} \!\!12
 +\sum_{n_2=0}^{\big\lfloor\frac{\la_1+2\la_2+\la_3-2}{2}\big\rfloor}
  \sum_{n_1=\big\lceil\frac{\la_1+n_2}{2}\big\rceil}^{\big\lfloor\frac{\la_1+n_2}{2}\big\rfloor}
  \sum_{n_3=\big\lceil\frac{\la_3+n_2}{2}\big\rceil}^{\big\lfloor\frac{\la_3+n_2}{2}\big\rfloor} \!\!6
 \nonumber\\[.2cm]
 &+\sum_{n_1=\max\big(\big\lceil\frac{\la_1}{2}\big\rceil,\big\lceil\frac{2\la_1+\la_2}{3}\big\rceil\big)}^{
   \big\lfloor\frac{3\la_1+2\la_2+\la_3-1}{4}\big\rfloor} \!\!4
 +\sum_{n_3=\max\big(\big\lceil\frac{\la_3}{2}\big\rceil,\big\lceil\frac{\la_2+2\la_3}{3}\big\rceil\big)}^{
  \big\lfloor\frac{\la_1+2\la_2+3\la_3-1}{4}\big\rfloor} \!\!4
 \nonumber\\[.2cm]
 &+\sum_{n_1=\big\lceil\frac{3\la_1+2\la_2+\la_3}{4}\big\rceil}^{\big\lfloor\frac{3\la_1+2\la_2+\la_3}{4}\big\rfloor}
  \sum_{n_2=\big\lceil\frac{\la_1+2\la_2+\la_3}{2}\big\rceil}^{\big\lfloor\frac{\la_1+2\la_2+\la_3}{2}\big\rfloor}
  \sum_{n_3=\big\lceil\frac{\la_1+2\la_2+3\la_3}{4}\big\rceil}^{\big\lfloor\frac{\la_1+2\la_2+3\la_3}{4}\big\rfloor}
  \!\!1,
\end{align}
which is seen to agree with (\ref{PA3}).

%

\end{document}